\documentclass[12pt]{article}

\usepackage[dvips]{graphicx}
\usepackage{makeidx}

\newfont{\pseudocode}{cmtt10}

\newtheorem{thm}{Theorem}
\newtheorem{cor}{Corollary}
\newtheorem{lem}{Lemma}

\newtheorem{assume}{Assumption}
\newtheorem{rmk}{Remark}

\newcommand{\dx}{{\dot x}}

\newcommand{\bfx}{\mathbf{x}}
\newcommand{\bfa}{\mathbf{a}}
\newcommand{\bfc}{\mathbf{c}}
\newcommand{\bfb}{\mathbf{b}}

\newcommand{\bfy}{\mathbf{y}}

\newcommand{\pd}{{\partial}}

\newcommand{\bfC}{\mathbf{C}}

\begin{document}

\baselineskip  7mm

\title{
\bf Parameter Estimation of Sigmoid Superpositions: Dynamical System Approach}
\author{Ivan Tyukin, Cees van Leeuwen, Danil Prokhorov $^{\$}$ \\ \\
\normalfont Laboratory for Perceptual Dynamics\\ RIKEN Brain Science Institute\\
2-1, Hirosawa, Wako-shi, Saitama, 351-0198,  Japan\\
tyukinivan@brain.riken.go.jp, ceesvl@brain.riken.go.jp
\\ \\ \normalfont $^\$$Ford Research Laboratory, \\ Dearborn, MI 48121, U.S.A.,
dprokhor@ford\@.com}
\date{}
\maketitle

\begin{abstract} Superposition of sigmoid function over a finite time
interval is shown to be equivalent to the linear combination of
the solutions of a linearly parameterized system of logistic
differential equations. Due to the linearity with respect to the
parameters of the system, it is possible to design an effective
procedure for parameter adjustment. Stability properties of this
procedure are analyzed.
\end{abstract}


\section{Introduction}

Static base functions are used in a variety of universal
function-approximation schemes. Their general form runs as
follows:  Let a given continuous function $g(t)$ be defined over a
compact time interval $[0,T]$. There will be a function $y(t)$,
represented as
\begin{equation}\label{approximation_furmula}
y(t)=\sum_{i=1}^n c_i f(a_i t + b_i),
\end{equation}
in which $f(\cdot): R\rightarrow R$ is a continuous function and
for any given $\varepsilon
> 0$, there are values of $n, a_i, b_i$, and $c_i$, such that for
all $t\in [0,T]$,
\[
|g(t)-y(t)|\leq \varepsilon.
\]
Among the functions $f(\cdot)$ for which approximation of $g(t)$
can be proven, the Gaussian and the sigmoid are the most
well-known ones. Approximation by sigmoid is often favored for,
amongst others, its very good rate of convergence with respect to
the number $n$ of additive terms in equation
(\ref{approximation_furmula}) \cite{Barron}. Recent results
\cite{Dingankar} have shown that
\[
\int_{0}^{t}(g(\tau)-y(\tau))^{2}d\tau = O
\left(\frac{1}{n^2}\right).
\]
Another advantage is that convergence is also possible in {\it
Sobolev space}, implying the existence of an optimal approximator
for derivatives of function $g(t)$ \cite{Hornik90};
\cite{Hornik94}.

In spite of significant progress in the fields of nonlinear
optimization and neural networks (a comprehensive review of a
neural learning algorithms is given in \cite{Haykin99}) an
estimation of the unknown values of parameters $a_i$, $b_i$, $c_i$
in (\ref{approximation_furmula}) is still a difficult problem.
Simple local optimization strategies, involving gradient descent,
fail to converge because of nonconvexity of the function with
respect to the parameters; global search algorithms
\cite{Kirkpatrick83}; \cite{Hansen92} are prohibitively expensive
computationally \cite{WalterProzanto97}, and second-order search
algorithms rely on assumptions relating to the error surface that
are not always met, for instance uniqueness of the extremum
\cite{Zangwill69}.

In order to address the parameter adjustment problem, simplifying
assumptions have been made \cite{Castillo02}. This approach, for
instance, requires that  the values of  each additive term $f(a_i
t + b_i)$ in (\ref{approximation_furmula}) over $[0,T]$ be known.
Under this assumption convergence to a global minimum could be
proven. The method was shown to have a very fast speed of
convergence. However, the requirement that the value of each term
be known imposes severe restrictions on the applicability of this
method. Following a different strategy, in recent years several
new methods have been proposed which are capable of avoiding local
minima by modifying the learning criterion (see, for instance
\cite{Lo2002}). Yet, these methods cannot guarantee that the
estimates of the unknown values of the parameters $a_i$, $b_i$,
$c_i$ converge to their true values (up to permutations). In our
view the underlying problem with these conventional methods is
that, whereas they use error minimization for approximating a
solution, they lack an explicit model of error dynamics. We will
propose a novel approach to estimate the values of the parameters
in (\ref{approximation_furmula}) utilizing elements of classical
control theory.

In this approach the values of function $g(t)$ are interpreted as
{\it reference signals}, the outputs of a dynamical system called
{\it reference system}. The reference signal is used in the
explicit definition of an error function as, for instance, the
difference with a {\it tracking signal}. This signal, in turn, is
considered the output $y(\theta,t): \Omega_{\theta}\times
R\rightarrow R$, $\theta\in\Omega_{\theta}$,
$\Omega_{\theta}\times R\rightarrow R$, $\theta\in\Omega_{\theta}$
of a dynamical system  called {\it tracking system} with parameter
vector $\theta=(\bfa^{T}:\bfb^{T}:\bfc^{T})$
 $\bfa,\bfb,\bfc \in R^{n}$ to be determined. Thus the problem
 of function approximation is transformed into
one of finding a suitable parameterization for a given tracking
system.

A similar strategy was used in \cite{Terekhov991},
\cite{Ampazis01} for different purposes. In these studies the
resulting equations remained nonlinear in their parameters. The
presently proposed transformation, however, will enable us to
represent the problem in terms of a nonlinear system that is
linear in its parameters. The linearity allows us to apply
conventional methods of adaptive control theory for stabilizing
the error dynamics and thus facilitate finding the optimal
solution. For this purpose, the learning problem is formulated as
one of adaptive tracking (or equivalently, synchronization between
reference and tracking system). To this problem we can apply the
method of Lyapunov functions, extending parameter space
$\Omega_{\theta}$ to $\{\alpha,\beta,\mathbf{C},\bfx(t)|
\alpha,\beta,\mathbf{C}\in R^{n}, \  \bfx(t): R\rightarrow R^{n}
\}$,
and use a simple rule for parameter adjustment in the enhanced
system dynamics. This provides us with a method potentially more
powerful than, for instance, gradient descent, which operates
entirely in the original parameter space by relying on the
contraction theorem.

It should be mentioned, however, that the problem of parameter
value identification has not completely been solved even for our
case of linearly parameterized, nonlinear systems. The solutions
available in the literature are formulated either for linear
systems \cite{Kreisselmeier77}; \cite{Magni95}; \cite{Basar96} or
for some special classes of nonlinear plants, assuming full state
measurement \cite{Didinsky95} or the possibility to transform the
system into an output injection form \cite{MarinoTomei93};
\cite{MarinoTomei95}. We do not wish to impose any such
restrictions. Instead we exploit the possibility to extend both
the reference and tracking signals to be repeated periodically
starting from the same initial conditions. By doing so we
significantly simplify the problem of searching for the optimal
values of unknown parameters.

A strategy similar to the one proposed is often used in iterative
learning control \cite{Arimoto84}; \cite{Arimoto88};
\cite{Messner91}; \cite{Phan98}  mostly for determining a
feed-forward control term which is defined as a function of time.
The time-variability of the solution severely reduces the
significance of these methods for our problem. Nevertheless, there
are several approaches that can be applied to search for unknown
parameters within an iterative learning control framework
\cite{Olerro90}; \cite{Hjalmarsson98}; \cite{Seo98}. These
approaches, however, according to our knowledge, are either
designed for linear dynamical systems or when dealing with
nonlinear systems cannot guarantee to stop at the non-local
solution. This motivates us not only to show the possibility to
transform the entire problem of static nonlinear optimization into
dynamics one but also to provide an algorithm to estimate the
unknown parameters of the resulting linearly parameterized system
of nonlinear differential equations.

The first step in our approach will be the selection of a ``base
function" for the reference and tracking systems, suitable for
representing a broad class of functions. We have chosen the
logistic differential equation 
\cite{Strogatz}. We will start off by providing an existence proof
for approximation in this system. The next step will be the
specification of an algorithm for parameter adjustment that
effectively finds the optimal solution in an interesting domain of
functions. We consider this problem for systems with unperturbed
conditions as well as with time-varying parameters. The former
constitutes a method for representing scalar functions in  one
variable, for instance time; the latter provides a method for
representing functions with multiple inputs. Finally, the
viability of the approach is demonstrated in examples comparing it
to gradient descent.

The paper is organized as follows. In Section 2 we formulate the
problem and introduce the class of systems to be analyzed. In
Section 3 we investigate the dynamic abilities of the system and
prove the approximation properties of the system. In Section 4 we
introduce the schemes to adjust the unknown parameters of the
system. In Section 5 we discuss multi-dimensional approximation
problems and show the possibility to utilize the same technique
for approximation of a system of nonlinear differential equations
with arbitrary smooth right-hand sides. Section 6 contains
simulation results for illustrative examples. Section 7 concludes
the paper.

\section{Problem Formulation}

Although  the  sigmoidal function approximation scheme has several
attractive features, the most important obstacle on the way to its
implementation remains the absence of an algorithm that guarantees
convergence to an optimal solution.  We suggest a strategy to turn
the problem of searching for the parameter values of the static
nonlinear parameterized map $f(\bfa,\bfb,\bfc,t)$, $\bfa,\bfb,\bfc
\in R^n$ into one of searching for linear parameter values of a
system of nonlinear differential equations:
\begin{equation}\label{trans}
\dot{\bfx}=\sum_{i=1}^n\xi_{1,i}(\bfx)\alpha_i + \sum_{i=1}^n
\xi_{2,i}(\bfx)\beta_i, \ y(\bfx)=\mathbf{C}\bfx,
\end{equation}
where $\bfx\in R^{n}$,
$\alpha=(\alpha_1,\dots,\alpha_n)^T,\beta=(\beta_1,\dots,\beta_n)^T\in
R^n$, $\xi_{1,i}: R^n\rightarrow R^n$, $\xi_{2,i}: R^n\rightarrow
R^n$ are continuous functions, $\mathbf{C}\in R^{n}$\footnote{We
would like to note that dimensions of the vectors $\alpha$ and
$\beta$ are not necessarily equal to $n$. Although we do not
discuss any other parameterization, a variety of alternative
descriptions with different parameterizations is possible.}.
Therefore, the first problem to be addressed is the existence of
such a transformation. The proposed solution uses  differential
logistic equations to realize system (\ref{trans}). This means we
will approach function $g(t)$ with a weighted sum $y(\bfx(t))$,
for which we then have to deal with the issue of identifying the
parameter values of (\ref{trans}). To this purpose, in
control-theoretic terms, system (\ref{trans}) is considered the
reference system, whereas the tracking system will have the
following description:
\begin{equation}\label{trans1}
\dot{\hat{\bfx}}=\sum_{i=1}^n\xi_{1,i}(\hat\bfx)\hat\alpha_i+\sum_{i=1}^n\xi_{2,i}(\hat\bfx)\hat\beta_i+\eta(y(\bfx),y(\hat{\bfx}),t),
\ y(\hat{\bfx})=\mathbf{\hat{C}}\hat{\bfx},
\end{equation}
where $\hat\bfx\in R^{n}$,
$\hat\alpha=(\hat{\alpha}_1,\dots,\hat{\alpha}_n)^T,\hat\beta=(\hat{\beta}_1,\dots,\hat{\beta}_n)^T\in
R^n$, $\mathbf{\hat{C}}\in R^{n}$.  Note the similarity in
structure between tracking and reference system, except for an
error function $\eta: R^3\rightarrow R^n$, added to the tracking
system. In what follows  symbols $\bfx(t)$, $\hat\bfx(t)$ denote
the solutions of differential equations (\ref{trans}),
(\ref{trans1}) with parameters $\alpha$, $\beta$ ($\hat\alpha$ and
$\hat\beta$) and starting from initial conditions $\bfx_0$. Sometimes
in order to stress this dependence explicitly we will write
$\bfx(\alpha,\beta,\bfx_0,t)$ or
$\hat\bfx(\hat\alpha,\hat\beta,\bfx_0,t)$.

 As both the reference and
tracking systems are described in the same manner, it is natural
to consider the combined system, which couples reference to
tracking system via output $y(\bfx(t))$ through the error function
$\eta(y(\bfx),y(\hat\bfx))$:
\begin{eqnarray}\label{ref_track}
\dot{\bfx}&=&\sum_{i=1}^n\xi_{1,i}(\bfx)\alpha_i+\sum_{i=1}^n\xi_{2,i}(\bfx)\beta_i,
\
y(\bfx)=\mathbf{C}\bfx, \nonumber \\
\dot{\hat{\bfx}}&=&\sum_{i=1}^n\xi_{1,i}(\hat\bfx)\hat\alpha_i+\sum_{i=1}^n\xi_{2,i}(\hat\bfx)\hat\beta_i+\eta(y(\bfx),y(\hat{\bfx}),t),
\ y(\hat{\bfx})=\mathbf{\hat{C}}\hat{\bfx}.
\end{eqnarray}
It is possible then to estimate the unknown parameters
$\alpha,\beta, \mathbf{\bfC}$ of the reference system. We start
out by assuming that the only uncertainties are in the vectors
$\alpha$ and $\beta$, while vector $\bfC$ is supposed to be known.
We will
propose an algorithm for parameter adjustment that is capable of
finding the solution. Our learning algorithm will belong to the
following class:
\begin{eqnarray}
\dot{\hat{\alpha}}&=&\mathcal{A}(y(\bfx),y(\hat{\bfx}),\bfx);
\nonumber \\
\dot{\hat{\beta}}&=&\mathcal{B}(y(\bfx),y(\hat{\bfx}),\bfx),
\end{eqnarray}
where operators $\mathcal{A}(\cdot)$ and $\mathcal{B}(\cdot)$ are
to be determined on the basis of the speed-gradient algorithm
\cite{Fradkov79}. If this strategy works, an extension would be to
consider cases where the reference system does not represent
function $g(t)$ completely (i.e. systems with unmodeled dynamics).

Thus, the questions to be addressed are: is it possible (at least
in theory) to transform a problem of nonlinear static optimization
into a problem of searching for linearly parameterized nonlinear
differential equations? If so, then how to estimate the parameters
of this nonlinear dynamical system in order to obtain qualitative
approximation? The next sections will provide us with the answers.

\section{Approximation with Logistic Differential Equations}

Let the following system be given:
\begin{eqnarray}\label{logistic_equations}
\dx_1&=&\alpha_1 x_1 (1- \beta_1 x_1); \nonumber \\
\dx_2&=&\alpha_2 x_2 (1- \beta_2 x_2); \nonumber \\
\cdots&=& \cdots \nonumber \\
\dx_n&=&\alpha_n x_n (1- \beta_n x_n); \nonumber \\
y({\bfx})&=&\mathbf{C}^T \bfx=\sum_{i}c_i x_i,  \ \
x_i(0)=\Delta_i,
\end{eqnarray}
where $\bfx=(x_1,\dots,x_n)^{T} \in R^n$ is a state vector,
$\alpha_i \in R$,  are parameters of system
(\ref{logistic_equations}), $y$ is an output function,
$\mathbf{C}=(c_1,\dots,c_n)^{T} \in R^n$ is a vector of parameters
associated with output $y$,  $x_i(0)\in R$ are initial conditions.

We begin our investigation by asking the question: what dynamics
can the autonomous system (\ref{logistic_equations}) produce as a
function of $t$? The answer to this question is formulated in the
following theorem:

\begin{thm}\label{logistic_approximation} Let continuously differentiable function  $g(t): R\rightarrow R$  be
given. Then for any $\varepsilon > 0$, $0<T<\infty$ and
$t\in[0,T]$ there are such numbers $n$, $\alpha_i$, $\beta_i$,
$c_i$ and initial conditions $x_i(0)=\Delta_i$ that the following
inequality holds:
\[
|y(\bfx(t))-g(t)|\leq \varepsilon.
\]
\end{thm}
Theorem \ref{logistic_approximation} proof is quite
straightforward and is based on the known fact that solution of
the logistic differential equation of the first order can be given
by a sigmoidal function \cite{Luenberger79}. Nevertheless, in
order to make the paper self-contained we present the proof in the
Appendix. Proofs of the subsequent theorems and lemmas are given
in the Appendix as well.

\begin{rmk}\label{remark1} \normalfont It follows from Theorem \ref{logistic_approximation}
proof that it is possible to transform the problem of nonlinear
function approximation by static sigmoidal functions into a
problem of choosing initial conditions and parameters $\alpha_i$
and $c_i$  of {\it dynamical} system (\ref{logistic_equations}),
where parameters ${\alpha}_i$  enter (\ref{logistic_equations})
linearly. One can observe, in addition, that under an appropriate
linear transformation $x_i\rightarrow x_i/c_i$ ($c_i\neq 0$) we
can get rid of uncertainties in $\mathbf{C}$ (see Remark
\ref{trans_rem_C} after Lemma \ref{lemma_transform} in Appendix 1)
and replace system (\ref{logistic_equations}) by
\begin{eqnarray}\label{logistic_equations1}
\dx_i&=&\alpha_i x_i +  \beta_i x^2_i; \nonumber \\
y({\bfx})&=&\sum_{i}x_i,  \ \ x_i(0)=\Delta_i/c_i,
\end{eqnarray}
where $\alpha_i$ and $\beta_i$ are to be determined. We formulate
this
\begin{cor}\label{logistic_approximation_cor} Let system (\ref{logistic_equations1}) and continuous differentiable  function  $g(t): R\rightarrow R$ be
given. Then for any $\varepsilon > 0$, $0<T<\infty$ and
$t\in[0,T]$ there are such numbers $n$, $\alpha_i$, $\beta_i$ and
initial conditions $x_i(0)$ that the following inequality holds:
\[
|y(\bfx(t))-g(t)|\leq \varepsilon.
\]
\end{cor}
This result will allow us to turn the problem of determining the
nonlinear parameters of a static function into a problem of
determining the linear parameters ${\alpha}_i$, $\beta_i$ of
system (\ref{logistic_equations1}). The restrictions are that the
values $x_i (0)$ will have to be known.
\end{rmk}

\begin{rmk}\label{remark31} \normalfont Theorem \ref{logistic_approximation}
proves that there is a one-to-one transformation of a function
approximation problem in terms of static sigmoidal functions to
one in terms of differential logistic equations. The latter,
therefore, shares all the advantages of the former, including the
very good convergence rate \cite{Dingankar} and its application in
Sobolev space \cite{Hornik94}.
\end{rmk}
Theorem \ref{logistic_approximation} merely states the existence
of parameters $\alpha_i$ and ${c}_i$ of system
(\ref{logistic_equations}) (or $\alpha_i$ and $\beta_i$ of system
(\ref{logistic_equations1})) that ensure arbitrarily small errors
between the system output and the reference function $g(t)$. It
does not answer the question how to derive the parameters.
However, the linearity of the system in its parameters simplifies
our task. We will show in Section 5 that in the multidimensional
case the resulting system will be linearly parameterized as well.
In the next section we will turn to the issue of how to find the
values of the parameters $\alpha_i$ that yield minimum errors.

\section{Parameter Adjustment Algorithm}

The question is whether it is possible to estimate the unknown
parameter values ${\alpha}_i$, $\beta_i$ for which
$g(t)-y(\hat{\bfx}(t))=0$ for $t\in[0,T]$, utilizing the linear
parameterization of system (\ref{logistic_equations1}). For
designing the estimation algorithm the following strategy was
used: first, it is assumed that the only uncertainties are in the
linear parameters $\alpha_i$, $\beta_i$, initial conditions
$\bfx(0)$ are assumed to be known. We formulate this in Assumption
\ref{assume_C}. First, our main algorithm is presented. Second,
after this algorithm is given we extend it to the cases where the
reference system does not represent the function $g(t)$
completely, i.e., with unmodeled dynamics. It will be possible to
invoke Theorem \ref{logistic_approximation} and show that any
function that merely is approached by reference system dynamics
can still effectively be modelled by the tracking system, albeit
within a margin of tolerance.

In order to proceed with the analysis we would like to introduce
the following assumption:
\begin{assume}\label{assume_C} Let continuous function $g(t)$, number of equations $n$ and initial
conditions $x_i(0)$ be given. There exist such parameter values
$\alpha_i$  and $\beta_i$ that for any $t\in[0,T]$ the following
equality holds for system (\ref{logistic_equations1}) solutions:
\[
g(t)-\sum_{i=1}^n c_i x_i(t)=0.
\]
\end{assume}
Assumption \ref{assume_C} states that the reference signal $g(t)$
can be represented by the output of system
(\ref{logistic_equations1}):
\[
g(t)=\sum_{i=1}^n c_i x_i(\alpha_i,\beta_i,x_i(0),t).
\]
The coefficients $c_i$ can be equal to the unity.

 In order to make the presentation more clear and compact,
 we would like to introduce a notational assumption
regarding the tracking and reference systems. Let us redefine the
system equations, denoting the right-hand side of
(\ref{logistic_equations1}) by
$\sum_{i=1}^n\xi_{1,i}(\bfx)\alpha_i +
\sum_{i=1}^n\xi_{2,i}(\bfx)\beta_i$, where
\[
\xi_{1,i}(\bfx)=\left(\begin{array}{l}
                    (i-1)\left\{
                          \begin{array}{l}
                           0\\
                           \cdots\\
                           0
                          \end{array}
                    \right.\\
                    x_i\\
                    (n-i)\left\{
                         \begin{array}{l}
                          0\\
                          \cdots\\
                          0
                         \end{array}
                    \right.
                \end{array},
                \right), \ \
\xi_{2,i}(\bfx)=\left(\begin{array}{l} (i-1)\left\{
                          \begin{array}{l}
                           0\\
                           \cdots\\
                           0
                          \end{array}
                    \right.\\
                    x^2_i\\
                    (n-i)\left\{
                         \begin{array}{l}
                          0\\
                          \cdots\\
                          0
                         \end{array}
                    \right.
                      \end{array}
                \right).
\]
Then both reference and tracking system can be rewritten in the
compact form (\ref{ref_track}) introduced in Section 2:
\begin{eqnarray}
\dot{\bfx}&=&\sum_{i=1}^n\xi_{1,i}(\bfx)\alpha_i+\sum_{i=1}^n\xi_{2,i}(\bfx)\beta,
\
y(\bfx)=\mathbf{C}\bfx, \nonumber \\
\dot{\hat{\bfx}}&=&\sum_{i=1}^n\xi_{1,i}(\hat\bfx)\hat\alpha_i+\sum_{i=1}^n\xi_{2,i}(\hat\bfx)\hat\beta_i+\eta(y(\bfx),y(\hat{\bfx}),t),
\ y(\hat{\bfx})=\mathbf{\hat{C}}\hat{\bfx}, \nonumber
\end{eqnarray}
where $\mathbf{C}=\hat\mathbf{C}=(1,\dots,1)^T$. Hence, to
complete the definitions of reference and tracking systems one
needs to determine $\eta(y(\bfx),y(\hat{\bfx}),t)$. One
possible way to do this is to define the function
$\eta(y(\bfx),y(\hat{\bfx}),t)$ as follows:
\[
\eta(y(\bfx),y(\hat{\bfx}),t)=K(t)(y(\hat\bfx)-y({\bfx})),
\]
where $K(t)=(k_1(t),\dots,k_n(t))^T$ and $k_i(t)$ are to be
specified later. The reason for such a structure is that we need
the tracking system ``to copy" the reference dynamics along a
manifold $y(\bfx)-y(\hat\bfx)=0$. Thus, an aggregated system which
contains both the reference system for signal $g(t)$ and tracking
system (\ref{logistic_equations1}) can be written in the following
form:
\begin{eqnarray}\label{ONO1}
\dot{\bfx}&=&\sum_{i=1}^n\xi_{1,i}(\bfx)\alpha_i+\sum_{i=1}^n\xi_{2,i}(\bfx)\beta_i,
\
y(\bfx)=\mathbf{C}\bfx, \nonumber \\
\dot{\hat{\bfx}}&=&\sum_{i=1}^n\xi_{1,i}(\hat\bfx)\hat\alpha_i+\sum_{i=1}^n\xi_{2,i}(\hat\bfx)\hat\beta_i+K(t)(y(\hat\bfx)-y({\bfx})),
\ y(\hat{\bfx})=\mathbf{\hat{C}}\hat{\bfx},
\end{eqnarray}


As has been mentioned in the beginning of the section, we would
like to obtain such estimates of the parameters $\alpha_i$,
$\beta_i$, that $g(t)-y(\hat{\bfx(t)})=0$ over time-interval
$[0,T]$. It was proposed in Section 2 to utilize conventional
speed-gradient like techniques to design the learning or
adaptation rule. For these methods, the parameters are supposed to
be adjusting on-line, that is in the same time-scale as the
reference and tracking systems evolve.  In general, it may take
much more time than $T$ (the length of the interval $[0,T]$) for
the estimates $\hat\alpha_i$, $\hat\beta_i$ to converge to
$\alpha_i$, $\beta_i$. However, the function $g(t)$ may not be
defined for $t>T$, and even if it is well defined over
$[T,\infty)$ then equivalence
$y(\bfx(\alpha,\beta,t,\bfx_0))\equiv
y(\hat{\bfx}(\hat\alpha,\hat\beta,t,\bfx_0))$ for $t>T$ does not
imply that $g(t)=y(\hat{\bfx}(\hat\alpha,\hat\beta,t,\bfx_0))$ for
any $t\in [0,T]$.

In addition, we note that logistic equations
(\ref{logistic_equations}) can be very unstable and may have
finite escape time depending on the vectors $\alpha$ and $\beta$.
For the reference system  this is not important as we assumed that
every solution $x_i$ of (\ref{logistic_equations}) can be
described by a sigmoid function and therefore is bounded. For the
tracking system, however, stability becomes very crucial. It is
very well possible that during $\hat\alpha$ and $\hat\beta$
adjustment and due to the term $K(t)(y(\bfx)-y(\hat{\bfx}))$ in
(\ref{ONO1}) the state $\hat{\bfx}$ of the reference system can
reach infinity in finite time thus making the whole system
unstable.

Taking these considerations into account, it is necessary to
redesign the reference and tracking systems in such a way that: 1)
$y(\bfx(\alpha,\beta,t,\bfx_0))\rightarrow
y(\hat{\bfx}(\hat\alpha,\hat\beta,t,\bfx_0))$ as
$t\rightarrow\infty$ implies that $|g(t)-
y(\hat{\bfx}(\hat\alpha,\hat\beta,t,\bfx_0))|<\varepsilon$ for any
$\varepsilon>0$ and arbitrary $t\in[0,T]$; and 2) the state
$\hat\bfx$ of the tracking system remains bounded for any $t>0$.

Our proposed solution to problem 1) is to let the reference signal
$g(t)$ be repeated periodically (see Fig. 1, where the initial
signal $g(t)$ is extended periodically along axis $t$).
Periodicity can be achieved by introducing special terms
($\lambda$ and $\sigma$ below) into the systems right-hand sides
that will periodically  force the states to move to $\bfx_0$ (with
period $T_1=T+\Delta T_2$, where $\Delta T_2$ is amount of time
needed to reach $\bfx_0$). In order to solve problem 2) we have to
make sure that state $\hat\bfx$ of the tracking system is bounded
for any $t>0$. This can be achieved if we force the states of {\it
both} systems to move to $\bfx_0$ as soon as $\|\hat\bfx\|$
exceeds certain bound $D$. Roughly speaking, one can add
time-varying negative feedback to both reference and tracking
systems, thus making the point $\bfx_0$  globally asymptotically
stable for both systems and, in addition, allowing the output
$y(\bfx(\alpha,\beta,t,\bfx_0))$ of the reference system to
coincide periodically with the segments of trajectory $g(t)$
defined over $[0,T]$.

In order to satisfy these requirements we introduce the next
\begin{assume}\label{assume_periodic} There is a positive
constant $l_0>0$ and function $\lambda: R^2\rightarrow R$
\[
\lambda (t, D)=\left\{ \begin{array}{ll}
                  0,& t  \in [(j-1) T_1, j T_1 - \Delta T_2 ) \ and \ \|\hat\bfx(t)\|<D \\
                  1,& t  \in [j T_1 - \Delta T_2, j T_1) \ or \ \|\hat\bfx(t)\|\geq D\
                  \end{array},  \ j\in\{1,2,\dots,\infty\},
            \right.
\]
such that the reference signal is given by the following system:
\begin{eqnarray}
& & \dot{\bfx}=\left(\sum_{i=1}^n\xi_{1,i}(\bfx)\alpha_i +
\sum_{i=1}^n\xi_{2,i}(\bfx)\beta_i\right)(1-\lambda(t,D))-
\lambda(t,D)\sigma(\bfx-\bfx(0))\nonumber \\
& & y(\bfx(t))=\tilde g(t),\nonumber
\end{eqnarray}
where $\sigma(\cdot)$ is a signum function:
\[
\sigma(\cdot)=(\sigma_1(\cdot),\dots,\sigma_n(\cdot))^T : \
\sigma_i(\bfx-\bfx(0))= \left\{
                 \begin{array}{ll}
                 1, & x_i-x_i(0)>0; \\
                 0, & x_i-x_i(0)=0; \\
                 -1 & x_i-x_i(0)<0,
                 \end{array}
           \right.
\]
$l_0\geq D/\Delta T_2$, $\tilde g(t_1)$, $t_1\in[0,\infty)$ is an
 extension of $g(t)$, $t\in[0,T]$ and $T_1=T+ \Delta T_2$.
\end{assume}

Assumption \ref{assume_periodic} requires an inclusion of several
extra parameters and functions into the generating system
right-hand side. Additional restrictions are to be introduced just
to make sure that for each $t=jT_1$, the following holds:
\[
x_i(t)=x_i(jT_1)=\hat{x}_i(t)=\hat{x}_i(jT_1)=x_i(0), \ \ j=\{1,2,\dots,\infty\}, \ \
i\in\{1,\dots,n\}.
\]

Taking into account Assumption \ref{assume_periodic} and the fact
that the tracking system is designed to copy the structure of the
reference system, we can write the combined reference and tracking
systems as follows:
\begin{eqnarray}\label{ONO2}
& & \dot{\bfx}=\left(\sum_{i=1}^n\xi_{1,i}(\bfx)\alpha_i +
\sum_{i=1}^n\xi_{2,i}(\bfx)\beta_i\right)(1-\lambda(t,D))-
\lambda(t,D)l_0 \sigma (\bfx-\bfx(0))  \nonumber  \\
& &
\dot{\hat\bfx}=\left(\sum_{i=1}^n\xi_{1,i}(\hat\bfx)\hat\alpha_i +
\sum_{i=1}^n\xi_{2,i}(\hat{\bfx})\hat\beta_i+K(t)(y(\hat\bfx)-y(\bfx))\right)(1-\lambda(t,D))-
\lambda(t,D)l_0 \sigma (\hat\bfx-\bfx(0))  \nonumber  \\
& & \hat y(t)= y (\hat{\bfx}(t))=\hat \bfC^T \hat{\bfx}(t);
\nonumber \\
& & y(t) = y (\bfx(t))=\bfC^T
 {\bfx}(t).
\end{eqnarray}
Before we introduce an adjustment rule for the tracking system let
us formulate the following lemma:

\begin{lem}\label{CEquiv} Let system (\ref{ONO2}) be given and $\hat{\mathbf{C}}^T\neq 0$. Consider
\[
|\hat\mathbf{C}^T\sum_{i=1}^{n}\left(\alpha_i(\xi_{1,i}(\hat{\bfx})-\xi_{1,i}(\bfx))+\beta_i(\xi_{2,i}(\hat{\bfx})-\xi_{2,i}(\bfx))\right)(1-\lambda(t,D))|+\epsilon\sum_{i=1}^n
k_i \hat{c}_i
\]
Then for any given constant $\delta>0$ there exist $k_i=k_i^\ast
\in R$  such that
\begin{eqnarray}\label{CEquiv_inequality}
|\hat\mathbf{C}^T\sum_{i=1}^{n}\left(\alpha_i(\xi_{1,i}(\hat{\bfx})-\xi_{1,i}(\bfx))+\beta_i(\xi_{2,i}(\hat{\bfx})-\xi_{2,i}(\bfx))\right)(1-\lambda(t,D))|+\epsilon\sum_{i=1}^n
k_i^{\ast}\hat{c}_i<0
\end{eqnarray}
for any $\epsilon >\delta$.
\end{lem}

According to Lemma \ref{CEquiv} for any positive $\delta>0$ the
existence of the coefficients $k^{\ast}_i$ satisfying inequality
(\ref{CEquiv_inequality}) is guaranteed. This property is very
important for the subsequent analysis. In fact, it states that the
error function $e=\hat{y}(t)-y(t)$
is attracted to the domain $|e|\leq \delta$ at
$\hat\alpha=\alpha$, $\hat\beta=\beta$, $\lambda(t,D)=0$ and
$k_i(t)=k^{\ast}_i$ as
\[
\dot e =
\left(\hat\mathbf{C}^T\sum_{i=1}^{n}\left(\alpha_i(\xi_{1,i}(\hat{\bfx})-\xi_{1,i}(\bfx))+\beta_i(\xi_{2,i}(\hat{\bfx})-\xi_{2,i}(\bfx))\right)+e\sum_{i=1}^n
k_i^{\ast}\hat{c}_i\right)(1-\lambda(t,D))
\]
and
\[
\frac{d}{dt}(0.5 e^{2})=e\dot{e}<0, \ \forall |e|>\delta.
\]

Let us introduce the adjustment rules for parameters
$\hat\alpha_i$, $\hat{\beta}_i$:
\begin{eqnarray}\label{learning_rule}
\dot{\hat{\alpha}}_i&=&-\gamma
e(t)S_\delta(e)\hat{\mathbf{C}}^T\xi_{1,i}(\hat\bfx)(1-\lambda(t,D)),\nonumber \\
\dot{\hat{\beta}}_i&=&-\gamma
e(t)S_\delta(e)\hat{\mathbf{C}}^T\xi_{2,i}(\hat\bfx)(1-\lambda(t,D)),
\\
S_\delta(e)&=&\left\{\begin{array}{ll}
                    1, & |e| > \delta\\
                    0, & |e|\leq \delta
                   \end{array}
            \right. \nonumber.
\end{eqnarray}
where $e(t)=\hat{y}(t)-y^{\ast}(t)$ is the tracking error,
$\gamma>0$ is a positive constant.

The stability properties of system (\ref{ONO2}) with algorithm
(\ref{learning_rule}) are formulated in:
\begin{thm}\label{learning_theorem} Let Assumptions
\ref{assume_C}, \ref{assume_periodic} hold, vector
$\hat{\mathbf{C}}\neq 0$, and function
$K(t)=(k_1(t),\dots,k_n(t))^T$ in (\ref{ONO2}) be given by the
following system of differential equations
\begin{eqnarray}\label{K_adjust}
\dot{k}_i=-\gamma S_\delta(e)e^2\hat{c}_i(1-\lambda(t,D)).
\end{eqnarray}
Then for any positive $\gamma>0$ all trajectories of system
(\ref{ONO2}) are bounded, and there exists $t_1>0$ such that for
any $t>t_1$ the following inequality holds:
\[
|y(\bfx)-y(\hat\bfx)|<\delta+\delta_1, \delta_1>0.
\]
\end{thm}

\begin{rmk}\label{parameter_estimation} \normalfont  Theorem \ref{learning_theorem} guarantees
that function $e(t)\lambda(t,D)$  in system (\ref{ONO2}) converges
to the domain $|e(t)\lambda(t,D)|<\delta$, where  constant
$\delta$ is defined in learning algorithm (\ref{learning_rule}).
Formally, $|e(t)\lambda(t,D)|<\delta$  does not automatically
imply that estimates $\hat\alpha$,$\hat{\beta})$ converge to the
point $\hat{\alpha}=\alpha$, $\hat\beta=\beta$ in the parameter
space. Nevertheless, according to formula (\ref{lt_2}) (see
Appendix, proof of Theorem \ref{learning_theorem}), one can derive
the following  estimate of how close we are to the solution
\begin{eqnarray}\label{par_est}
& & \|\hat{\alpha}(t_0)-\alpha\|^{2}_{\gamma^{-1}} +
\|\hat{\beta}(t_0)-\beta\|^{2}_{\gamma^{-1}}  -
\|\hat{\alpha}(t)-\alpha\|^{2}_{\gamma^{-1}} -
\|\hat{\beta}(t)-\beta\|^{2}_{\gamma^{-1}}
\nonumber\\
&
&\geq\|K(t)-k^{\ast}\|^{2}_{\gamma^{-1}}-\|K(t_0)-k^{\ast}\|^{2}_{\gamma^{-1}}+
2\int_{t_0}^{t}S_{\delta}(e)|e(\tau)\lambda(\tau,D)\sum_{j=1}^{n}{k^{\ast}_j}\hat{c}_j|\delta_1
d\tau.
\end{eqnarray}
Equation (\ref{par_est}) may be taken to reflect the quality of
estimation of the unknown parameters $\alpha$ and $\beta$. In
particular, if we choose $K(t_0)=0$, then
\begin{eqnarray}
& & \|\hat{\alpha}(t_0)-\alpha\|^{2}_{\gamma^{-1}} +
\|\hat{\beta}(t_0)-\beta\|^{2}_{\gamma^{-1}}  -
\|\hat{\alpha}(t)-\alpha\|^{2}_{\gamma^{-1}} -
\|\hat{\beta}(t)-\beta\|^{2}_{\gamma^{-1}}
\nonumber\\
&
&\geq\|K(t)-k^{\ast}\|^{2}_{\gamma^{-1}}-\|k^{\ast}\|^{2}_{\gamma^{-1}}+
2\int_{t_0}^{t}S_{\delta}(e)|e(\tau)\lambda(\tau,D)\sum_{j=1}^{n}{k^{\ast}_j}\hat{c}_j|\delta_1
d\tau.\nonumber
\end{eqnarray}
Therefore, the smaller the norm $\|K(t)\|$, the greater is the
chance that the difference
\begin{eqnarray}\label{dev_pos}
& & \|\hat{\alpha}(t_0)-\alpha\|^{2}_{\gamma^{-1}} +
\|\hat{\beta}(t_0)-\beta\|^{2}_{\gamma^{-1}}  -
\|\hat{\alpha}(t)-\alpha\|^{2}_{\gamma^{-1}} -
\|\hat{\beta}(t)-\beta\|^{2}_{\gamma^{-1}}.
\end{eqnarray}
is nonnegative. On the other hand, given the values of
$\delta$, $\delta_1$, $D$,
 $\hat{\mathbf{C}}$ and bounds for $\alpha$, $\beta$, one can explicitly estimate vector
 $k^{\ast}$,
satisfying inequality (\ref{CEquiv_inequality}). Hence in this
case formula (\ref{par_est}) gives explicit bounds for the
deviations of the estimates $\hat\alpha$, $\hat\beta$ with respect
to $\alpha$ and $\beta$. Furthermore, for known $k^{\ast}$ it is
possible to get rid of time-varying coefficients $k_i(t)$ in
(\ref{ONO2}), replacing them by $k^{\ast}_i$. In this case
difference (\ref{dev_pos}) is positive if $|e(t)\lambda(t,D)|$
exceeds $\delta$ at some time $t_1$.

In general in order to ensure the positiveness of difference
(\ref{dev_pos}) for a given parameterization of the reference
system, it is necessary to consider more carefully the dynamics of
the following deviation $\rho=\hat{\bfx}-{\bfx}$ at
$\hat\alpha=\alpha$ and $\hat\beta=\beta$ over the time intervals
where $\lambda(t,D)=0$:
\[
\dot{\rho}=\left(\sum_{i=1}^{n}\alpha_i(\xi_{1,i}(\hat{\bfx})-\xi_{1,i}(\bfx))+\beta_i(\xi_{2,i}(\hat{\bfx})-\xi_{2,i}(\bfx))\right)+K(t)\hat\mathbf{C}^{T}(\hat{\bfx}-\bfx)
\]
Functions $\xi_{1,i}$ and $\xi_{2,i}$ are differentiable with
respect to their arguments. Therefore there exist such
$\Xi_{1,i}(\hat\bfx,\bfx)$ and $\Xi_{2,i}(\hat\bfx,\bfx)$ that the
following equalities hold:
\begin{eqnarray}
& & \Xi_{1,i}(\hat\bfx,\bfx)(\hat\bfx-\bfx)=\xi_{1,i}(\hat{\bfx})-\xi_{1,i}(\bfx)\nonumber\\
& &
\Xi_{2,i}(\hat\bfx,\bfx)(\hat\bfx-\bfx)=\xi_{2,i}(\hat{\bfx})-\xi_{2,i}(\bfx)\nonumber
\end{eqnarray}
Then derivative $\dot\rho$ can be written in the following form
\begin{eqnarray}\label{discuss_1}
\dot{\rho}=\left(\sum_{i=1}^{n}\alpha_i\Xi_{1,i}(\hat\bfx,\bfx)+\beta_i\Xi_{2,i}(\hat\bfx,\bfx)
+K(t)\hat\mathbf{C}^{T}\right)\rho
\end{eqnarray}
It can be derived  from Theorem \ref{learning_theorem} proof that
the existence of a positive function $V(y(\hat\bfx),y(\bfx))$ with
time derivative $\dot{V}$ at $\hat\alpha=\alpha$,
$\hat\beta=\beta$ satisfying,
\begin{eqnarray}\label{proof_ineq}
\dot{V}(y(\hat\bfx),y(\bfx))=-W(y(\hat\bfx)-y(\bfx))
\end{eqnarray}
where $W(\cdot)$ is a positive definite function, guarantees
monotonic increase of the difference (\ref{dev_pos}). Therefore,
if one can find vector $K(t)$ such that it asymptotically
stabilizes system (\ref{discuss_1}) for the given domain of
parameters $\alpha$, $\beta$, and furthermore, inequality
(\ref{proof_ineq}) holds, then the positiveness of difference
(\ref{proof_ineq}) is guaranteed. The problem of determining
$K(t)$ however is not very easy to solve, especially for nonlinear
systems. Even for linear ones, a similar problem known in the
literature as the Brockett problem\footnote{Let the following
triplet of matrixes be given $A$, $B$, $C$ $\in$ ${R^{n\times
n}}$. Under what conditions does a time-variant matrix $K(t)$
exist such that system
\[
\dot{\bfx}=A\bfx+BK(t)C\bfx, \ x\in R^n
\]
is asymptotically stable?} \cite{Brockett98} has positive
solutions at present for  systems of second and third order
\cite{Leonov2000,Moreau99}. Nevertheless, despite the obvious
difficulties, we believe that the question of searching for the
suitable $K(t)$ ensuring inequality (\ref{proof_ineq}) for system
(\ref{discuss_1}) could be an achievable goal for future studies.
\end{rmk}

It is desirable to note that Theorem \ref{learning_theorem}
requires the validity of Assumption \ref{assume_C}. Assumption
\ref{assume_C} allowed us to model the function $g(t)$ by a
reference system of the same structure as the tracking one. This
feature
has been exploited in the proof of the theorem and played an
important role in order to guarantee convergence of errors to a
neighborhood of the origin. This assumption may be too restrictive
as it requires strict equivalence between reference and tracking
signals for $\hat\alpha=\alpha$, $\hat\beta=\beta$. We are now
ready to abandon this assumption by invoking Theorem
\ref{logistic_approximation} again.

If Assumption \ref{assume_C} does not hold this leads to nonzero
error  $\varepsilon (t)$ between the output
$y(\bfx)=\mathbf{C}^{T}\bfx(t)$ of the reference system
(\ref{ref_syst}) and signal $g(t)$ to be tracked:
\[
\varepsilon(t)=\sum_{i=1}^n {c}_i{x}_i(t) -g(t).
\]
Let us assume that $g(t)$ is continuously differentiable, then
$\varepsilon(t)$ is differentiable as well. We denote its first
derivative by $d\varepsilon(t)$:
\begin{equation}\label{logistic_unmodeled_error}
\frac{d}{dt}(y(\bfx(t))-g(t))=\sum_{i=1}^{n} {c}_i\dot{x}_i -
\dot{g}(t) =d\varepsilon(t).
\end{equation}
Due to the compactness of the interval $[0,T]$ we can conclude
that derivative $d\varepsilon(t)$ is bounded:
\[
|d\varepsilon(t)| < s.
\]

Let us derive the error
$e(t)=y(\hat\bfx)-g(t)=y(\hat\bfx)+\varepsilon(t)-y(\bfx)$
dynamics taking into account that $\mathbf{C}=\hat\mathbf{C}$ and,
in addition, that function $\varepsilon(t)$ can be considered as
an unmeasured disturbance subtracted from  the output $y(\bfx(t))$
generated by the reference system (\ref{ref_syst}):
\begin{eqnarray}\label{error_system2}
\dot{e}&=&\hat\mathbf{C}^{T}\left(\sum_{i=1}^{n}\hat\alpha_i\xi_{1,i}(\hat\bfx)-\alpha_i\xi_{1,i}(\bfx)+\hat\beta_i\xi_{2,i}(\hat\bfx)-\beta_i\xi_{2,i}(\bfx)\right)(1-\lambda(t,D))-d\varepsilon(t)\nonumber\\
& & +
\hat\mathbf{C}^{T}\left(K(t)(y(\hat\bfx)-y(\bfx)+\varepsilon(t))(1-\lambda(t,D))+l_0(\sigma(\bfx-\bfx_0)-\sigma(\hat\bfx-\bfx_0))\lambda(t,D)\right)
\end{eqnarray}
The only difference between error dynamics according to Assumption
\ref{assume_C}  and the expression given in (\ref{error_system2})
is in the term $d\varepsilon(t)+\mathcal{C}^{T}K(t)\varepsilon(t)$
which represents the unmodeled dynamics of $g(t)$.

There are several ways to deal with such an uncertainty. One of
them is to include a {\it dead-zone} into the parameter adjustment
scheme \cite{APoznyak} and chose $K(t)=const$. The algorithms with
a dead-zone will have the same form as (\ref{learning_rule}):
\begin{eqnarray}\label{learning_rule1}
\dot{\hat{\alpha}}_i&=&-\gamma
e(t)S_\delta(e)\hat{\mathbf{C}}^T\xi_{1,i}(\hat\bfx)(1-\lambda(t,D)),\nonumber \\
\dot{\hat{\beta}}_i&=&-\gamma
e(t)S_\delta(e)\hat{\mathbf{C}}^T\xi_{2,i}(\hat\bfx)(1-\lambda(t,D)),
\\
S_\delta(e)&=&\left\{\begin{array}{ll}
                    1, & |e| > \delta\\
                    0, & |e|\leq \delta
                   \end{array}
            \right. \nonumber.
\end{eqnarray}
except that the width $\delta$ of the dead-zone is to depend on
the bounds for $d\varepsilon(t)$ and
$\mathcal{C}^{T}K\epsilon(t)$. Theoretical analysis of the
stability of the whole system with learning rule
(\ref{learning_rule1}) can be done in the same manner as with
(\ref{learning_rule}).

It is clear that the tolerance of the resulting learning process
will depend on the dead-zone width $\delta$, which is exactly the
upper bound of $d\varepsilon(t)+\mathcal{C}^{T}K\epsilon(t)$.
Therefore, in general, applicability of the proposed learning
rules strongly depends on a smoothness of $\varepsilon(t)$ (in the
sense of the maximum absolute value of its first derivative). We
may deal with this issue by referring to the properties of this
approximation scheme in Sobolev space
\cite{Hornik94};\cite{Hornik90}. It can be shown that for any
arbitrary small $\delta_2>0$ there exists a network that can
approximate a given reference function $g(t)$ such that both
derivative $d\varepsilon(t)$ and $\varepsilon(t)$ satisfy the
following estimation:
$|d\varepsilon(t)+\mathcal{C}^{T}K\epsilon(t)|<\delta_2$. Hence,
learning algorithm (\ref{learning_rule1}) will still be applicable
even in the presence of nonzero differentiable error
$\varepsilon(t)$ between the reference signal and outputs of the
tracking system at $\hat{\alpha}=\alpha$, $\hat{\beta}=\beta$.
What value of $\delta$ is admissible will depend on the dimension
of the system.

\section{Discussion}{

Here we discuss multi-dimensional extensions with an eye for
possible neural network applications of our approach. Theorem
\ref{logistic_approximation} states that any continuous function
of $t$ can be approximated over time interval $[0,T]$ by a linear
combination of the solutions of system
(\ref{logistic_equations1}). It is desirable to note that we can
choose function $g(t)$ in such a way that the following equality
holds:
\begin{equation}\label{transform1}
g(t)=\tilde{g}(\xi(t)),
\end{equation}
where $g\in C^1$, $\xi(t)$ is a smooth function of $t$. Let us
suppose that system (\ref{logistic_equations1}) realizes function
$\tilde{g}(\xi)$. This means that
\[
\tilde{g}(\xi)=\sum_{i=1}^{n}{c}_i{x}_i(\xi),
\]
where $\dot{{x}}_i= {\alpha}_i {x}_i(1-\beta_i {x}_i)$. Then we
consider function $\tilde{g}(\xi)$  as a function of time $t$
which satisfies equation (\ref{transform1}). Therefore due to
formula (\ref{transform1}) we can write:
\[
\tilde{g}(\xi(t))=\sum_{i=1}^{n}{c}_i{x}_i(\xi(t)).
\]
Moreover
\[
\dot{g}(t)=\frac{d}{dt}\tilde{g}(\xi(t))=\frac{\pd}{\pd\xi}\tilde{g}(\xi)\frac{\pd}{\pd
t } \xi(t)=\sum_{i=1}^{n}{c}_i {x}_i(1-\beta_i{x}_i)\dot{\xi}.
\]
Hence under the following assumptions:
$\dot{g}(t)=\dot{\tilde{g}}(t)$ at $t=0$ and
$g(0)=\tilde{g}(\xi(0))$ we can see that linear combination
$\sum_{i=1}^n {c}_i {x}_i(t)$ of the solutions of system
\[
\dot{{x}}_i= {\alpha}_i {x}_i(1-\beta_i {x}_i)\dot{\xi}(t)
\]
realizes function $g(t)$ and vice-versa. This simple observation
suggests how to extend the result to the multi-dimensional case.
It is possible to consider a reference function
$g(\xi_1,\dots,\xi_m)$ with $m$ inputs as a function of time $t$:
$g(\xi_1(t),\dots,\xi_m(t))$. Then a system which realizes
function $g(\xi_1(t),\dots,\xi_m(t))$ can be represented in the
following form:
\begin{eqnarray}\label{logistic_multidimensional}
\dot{{x}}_i&=&\left(\sum_{j=1}^m\alpha_{i,j}\dot{\xi}_j(t)\right){x}_i(1-\beta_{i,j}{x}_i);
\nonumber \\
y(\hat{\bfx}(t))&=&\sum_{i=1}^n{c}_i{x}_i(t).
\end{eqnarray}
If we return to the approximation problem we  may observe on
account of Theorem \ref{logistic_approximation}  that system
(\ref{logistic_multidimensional}) is able to approximate a given
function $g(\xi_1,\dots,\xi_m)$ over a given compact domain in
such a way that for a particular trajectory
$(\xi_1(t),\dots,\xi_m(t))$ and any given constant $\varepsilon>0$
there exist parameters $\alpha_{i,j}$, $\beta_{i,j}$, ${c}_i$,
initial conditions and number $n$ satisfying the following:
\[
|g(\xi_1(t),\dots,\xi_m(t))-y({\bfx}(t))|\leq \varepsilon.
\]
Curve $\xi(t)$ should be designed in such a way that good
approximation along the curve $\xi(t)$ implies good approximation
along the whole surface. Intuitively, this depends on the degree
to which the curve "covers" the space. In other words, the more
complex curve $(\xi_1(t),\dots,\xi_m(t))$ is, the better the
approximation that can be achieved over the given compact
interval.

An important consequence of this description is that a system of
{\it coupled} logistic differential equations
(\ref{logistic_multidimensional}) may realize an approximation of
a nonlinear time-invariant system of the following type:
\begin{equation}\label{plant_approximation}
\dot{\bfy}=\chi(\bfy),
\end{equation}
where $\chi(\cdot): R^n\rightarrow R^n$ is an arbitrary smooth
function. Let us explain this. Denote:
\[
\mathcal{F}(\bfx,\bfb,\bfc,t)=\sum_{i=1}^{n}c_i f(a_i t +b_i).
\]
Consider system (\ref{logistic_multidimensional}) for $m=1$ and
replace $\dot{\xi}(t)$ by $\xi(t)$:
\begin{eqnarray}\label{logistic_feedback1}
\dot{{x}}_i&=&\alpha_{i}\xi(t){x}_i(1-\beta_i{x}_i);
\nonumber \\
y({\bfx}(t))&=&\sum_{i=1}^n{c}_i{x}_i(t)=\mathcal{F}(\alpha,\beta,\bfx_0,{\bfC},\int_{0}^t\xi(\tau)d\tau).
\end{eqnarray}
One may substitute function $y(t)$ in (\ref{logistic_feedback1})
instead of $\xi(t)$. This leads immediately to the following
equations:
\begin{eqnarray}\label{logistic_feedback2}
\dot{\hat{x}}_i&=&\alpha_{i}y(t){x}_i(1-\beta_i {x}_i);
\nonumber \\
y(t)&=&\mathcal{F}(\alpha,\beta,\bfx_0,{\bfC},\int_{0}^t
y(\tau)d\tau).
\end{eqnarray}
Denoting $z(t)=\int_{0}^t y(\tau)d\tau$ and taking into account
that $y=\sum_{i=1}^n{c}_i{x}_i$ we can rewrite system
(\ref{logistic_feedback2}) in the following manner:
\begin{eqnarray}\label{logistic_feedback3}
\dot{\hat{x}}_i&=&\alpha_{i} \left(\sum_{j=1}^n{c}_j{x}_j\right)
{x}_i(1-\beta_i{x}_i);
\nonumber \\
\dot{z}&=&\sum_{i=1}^n{c}_i{x}_i,
\end{eqnarray}
where the new output function $z(t)$ satisfies the following
differential equation:
\[
\dot{z}=\mathcal{F}(\alpha,\beta,\bfx_0,{\bfC},z).
\]
$\mathcal{F}(\alpha,\beta,\bfx_0,{\bfC},z)$ may realize function
$\chi(z)$ with given tolerance subject to the choice of the
parameters $\alpha,\beta,\bfx_0,{\bfC}$ and the number of
equations in (\ref{logistic_feedback3}). In the same fashion one
can derive the results for $m>1$ and obtain the corresponding
systems for differential equations:
\[
\dot{z}_i=\mathcal{F}_i(\alpha,\beta,\bfx_0,{\bfC},z_1,z_2,\dots,z_i,\dots,z_n),
\]
thus approximating (\ref{plant_approximation}).

There are two important observations to be made regarding system
(\ref{logistic_feedback3}). First, one may notice that system
(\ref{logistic_feedback3}) is a specific instance of the
Cohen-Grossberg model \cite{Cohen83}. Therefore, it is possible to
claim  that Cohen-Grossberg models of several differential
equations, each of which has relatively simple description (for
instance, coupled logistic differential equations), in principle,
are capable of approximating every nonlinear dynamical system with
smooth right-hand sides (subject to appropriate choice of the
number of differential equations, initial conditions and
parameters). Furthermore, the learning algorithms, introduced in
the paper can be applied to these models as well, and their
stability may be proven in the same fashion. Second, it is
desirable to notice that this approach allows us to introduce an
alternative learning technique to that of backpropagation through
time \cite{Werbos90}, albeit for continuous-time systems. A
detailed discussion of these topics is beyond the scope of the
present paper.

The algorithms introduced in the paper guarantee that under
certain circumstances the estimates $\hat{\alpha}$, $\hat{\beta}$
approach to a domain around $\alpha$, $\beta$. Still, they cannot
guarantee that $\hat{\alpha}\rightarrow\alpha$ and
$\hat\beta\rightarrow\beta$.  An interesting problem, therefore,
is  whether it is possible to design a tracking system that
guarantees convergence of $\hat{\alpha}$, $\hat{\beta}$ to
$\alpha$ and $\beta$ respectively. This problem in our opinion is
closely related to the problem of {\it adaptive observer design}
\cite{Marino95} for the reference system in (\ref{ONO2}):
\begin{eqnarray}\label{ref_syst}
& & \dot{\bfx}=\left(\sum_{i=1}^n\xi_{1,i}(\bfx)\alpha_i +
\sum_{i=1}^n\xi_{2,i}(\bfx)\beta_i\right)(1-\lambda(t,D))-
\lambda(t,D)\sigma(\bfx-\bfx(0))\nonumber \\
& & y(\bfx(t))=\mathbf{C}^T\bfx.
\end{eqnarray}
A prerequisite for applying the corresponding method is that these
systems are transformed into the {\it canonical observable form}
\cite{Bastin88}. For nonlinear systems that are linear in
parameters necessary and sufficient conditions for this have been
given  \cite{Marino90}. These conditions do not hold, however, for
the parameterizations of type (\ref{ref_syst}). Therefore, the
question remains open, whether is it possible to find such
linearly parameterized nonlinear system and corresponding output
function $y(\bfx)$, such that 1) its parameters can be transformed
by one-to-one mapping into those of sigmoid superposition, and 2)
the parameterization of this system obeys assumptions introduced
in work \cite{Marino90} (see Theorem 3.1). If one finds such a
suitable parameterization, then the problem of finding the ``true"
parameters (subject to permutations) can be solved  effectively.


\section{Examples}

In this section we illustrate the theoretical results with
examples. First we consider application of Theorem
\ref{learning_theorem} to the search for unknown parameter values
of a single sigmoid function and then  show the effectiveness of
our method in comparison with the conventional schemes for
two-dimensional optimization problem. In addition we illustrate
our method with the results of computer simulations performed for
a system consisting of 10 sigmoidal functions.

\subsection{Example 1}

Let us illustrate the possibility to search for the parameters
${\alpha}_i$ and ${c}_i$ simultaneously. As has been suggested in
Section 3, instead of the parameters $\alpha_i$ and ${c}_i$ we
will deal with ${\alpha}_i$ and ${\beta}_i=\alpha_i/{c}_i$.
Reference function $g(t)$ has been chosen to satisfy:
\[
g(t,\alpha,c)=\frac{c}{1+e^{-\alpha t +2.944}},
\]
where $c=2$, $\alpha=2/3$. We design the reference and tracking
systems as follows:
\begin{eqnarray}\label{example3}
& &\dot{x}=(\alpha x -{\beta}
{x}^2)(1-\lambda(t))-\lambda(t)(l_0\sigma(x-x(0))) \nonumber \\
& & \dot{\hat{x}}=(\hat\alpha \hat{x} -\hat{\beta}
\hat{x}^2)(1-\lambda(t))-\lambda(t)(l_0\sigma(\hat{x}-x(0)))-K(t)e,
\end{eqnarray}
where $\alpha=2/3$, ${\beta}=1/3$, $l_0=1$, $x(0)=0.1$,
$K(t)=0.2$, $e=\hat{x}-x$. Function $\lambda(t)$ was chosen to be
a periodic function with period $T=10$ sec, pulse width is $1$ sec
and unit amplitude (one may easily check that this parameter
setting ensures exact matching between function $g(t)$ and $x(t)$
over time interval $[0,9]$).

Adaptation rules to adjust the parameters $\hat{\alpha}$ and
$\hat{\beta}$ may be written as follows:
\begin{eqnarray}\label{example3_learning}
\dot{\hat{\alpha}} & = & -0.2 e(t)\hat{x}(t)(1-\lambda(t));
\nonumber \\
\dot{\hat{\beta}} & = & 0.2 e(t)\hat{x}^2(t)(1-\lambda(t)).
\end{eqnarray}
In order to make the example more illustrative we would like to
compare the performance of algorithm (\ref{example3_learning})
with a conventional pattern-by pattern gradient scheme:
\begin{eqnarray}\label{example3_gradient}
\dot{\hat{\alpha}}&=&-0.2 e(t)\frac{\pd
g(t,\hat{\alpha},\hat{c})}{\pd
\hat{\alpha}} \nonumber \\
\dot{\hat{c}}&=&-0.2 e(t)\frac{\pd g(t,\hat{\alpha},\hat{c})}{\pd
\hat{c}}
\end{eqnarray}
and batch rule:
\begin{eqnarray}\label{example3_batch}
\dot{\hat{\alpha}}&=&-0.2\frac{\pd J(\hat{\alpha},\hat{c})}{\pd
\hat{\alpha}} \nonumber \\
\dot{\hat{c}}&=&-0.2\frac{\pd J(\hat{\alpha},\hat{c})}{\pd
\hat{c}},
\end{eqnarray}
where
\[
J(\hat{\alpha},\hat{c})=\int_{0}^{9}(g(\tau,\hat{\alpha},\hat{c})-g(\tau,{\alpha}^{\ast},{c}^{\ast}))^2
d\tau
\]
Results of such a comparison are shown if Figures 2-5. In Figure 2
there are two trajectories of the parameters $\hat{\alpha}(t)$ and
$\hat{c}(t)$ in two-dimensional space. The first curve is obtained
from the trajectories of $\hat{\alpha}(t)=\hat{\alpha}(t)$,
$\hat{c}(t)=\hat{\alpha}(t)/\hat{\beta}(t)$ and results from
algorithm (\ref{example3_learning}) with initial conditions
$\hat\alpha(0)=-3, \ \hat{\beta}(0)=1$. Curve 2 is a solution of
(\ref{example3_gradient}) starting from initial conditions
$\hat\alpha(0)=-3$, $\hat{c}(0)=-3$. It can be seen that algorithm
(\ref{example3_learning}) reaches the global minimum. Conventional
gradient descent fails to do so. It appears unstable and goes
through a neighborhood of the global minimum along a valley. This
process is shown in Fig. 2. In addition, algorithm
(\ref{example3_learning}) is much faster than
(\ref{example3_gradient}) (see Fig. 3 for details).

Figure 4 reflects another interesting feature of algorithm
(\ref{example3_learning}). Whereas the conventional gradient
algorithm starting from $\hat{\alpha}(0)=3, \hat{c}(0)=-3$ goes
towards the goal along the isolines (Curve 2), algorithm
(\ref{example3_learning}) does not stick to isolines. Instead, it
goes through infinity in the coordinates $\hat{\alpha},\hat{c}$.
This is not because of any singularities with respect to the
coordinates $\hat{\alpha}$, $\hat{\beta}$ but is due simply to the
transformation $\hat{c}=\hat{\alpha}/\hat{\beta}$, when
$\hat\beta$ goes through zero.

Figure 5 contains the trajectories of the solutions obtained with
algorithm (\ref{example3_batch}). Curve 1 shows the trajectory
corresponding to initial conditions $\hat{\alpha}(0)=-3,
\hat{c}(0)=-3$, Curve 2 is related to initial conditions
$\hat{\alpha}(0)=3, \hat{c}(0)=-3$. It is easy to see that this
algorithm gets stuck in local minima.

The performance of algorithm (\ref{example3_learning}) is not
surprising because it uses information about the system properties
in a more intelligent way than gradient descent methods do. In
addition some coordinate transformation has been used and the
process of searching for the minimum is organized in a different
coordinate system. All the results relating to stability, however,
remain true for the functions which may be represented by a
superposition of sigmoid function only.

\subsection{Example 2}

In addition to the simple example of the previous section which
merely illustrates the design procedure for the parameters
adjustment rules proposed in the paper, we would like to present
more supporting results of computer simulation of our algorithms
for a larger number of functions in superposition. We consider the
sum of 10 sigmoid functions
\[
g(t,\alpha,\mathbf{C})=\sum_{i=1}^{10}\frac{c_i}{1+e^{-\alpha_i t
+b_i}},
\]
where parameters $b_i$ and $c_i$ are assumed to be known and $t\in
[0,T]$. According to the results presented, this sum is equivalent
to the solutions of the corresponding system of logistic equations
(\ref{logistic_equations}) with known $\beta_i$, $c_i$ and initial
conditions. The only uncertainties are in parameters $\alpha_i$.
First, we extend the reference signal $g(t)$ to be periodically
repeated over $[0,\infty)$:
\[
\tilde{g}(t)=\left\{\begin{array}{ll}
                     g(t), & t\leq T\\
                     0, & T<t<T+\Delta T_2, \\
                     g(t-T-\Delta T_2), & t> T+\Delta T_2
                    \end{array}\right.
\]
Then we design the tracking system
\begin{eqnarray}\label{track_ex}
\dot{\hat{x}}_i=\hat{\alpha}_i\hat{x}_i(1-\hat{x}_i)(1-\lambda(t,D))
+ k_i(t)e (1-\lambda(t,D))-\lambda(t,D)l_0\sigma(\hat{x}_i-x_i(0))
\end{eqnarray}
and adaptation algorithm
\begin{eqnarray}\label{alg_ex}
\dot{\hat\alpha}_i=-\gamma S_{\delta}(e)e
\hat{x}_i(1-\hat{x}_i)(1-\lambda(t,D))\nonumber \\
\dot{k}_i(t)=-\gamma S_{\delta}(e) e^2 c_i (1-\lambda(t,D))
\end{eqnarray}
where $D=10$ (taking into account that $|x_i|\leq 1$ we have to
choose $D>1$), $\lambda(t,D)$ is a $T+\Delta T_2$ periodic
function with the pulse width $\Delta T_2$, $\delta=0.0001$,
$\gamma=0.001$, $T=2$, $\Delta T_2=1$, $l_0=10$. Initial
conditions $x_i(0)$ and parameters $c_i$ were randomly chosen and
their exact values are given below:
\[
 \begin{array}{ll}
 x_1(0) & =0.1\\
 x_2(0) & =0.2\\
 x_3(0) & =0.3\\
 x_4(0) & =0.2\\
 x_5(0) & =0.5\\
 x_6(0) & =0.1\\
 x_7(0) & =0.7\\
 x_8(0) & =0.2\\
 x_9(0) & =0.6\\
 x_{10}(0) & =0.4
 \end{array}, \ \
 \begin{array}{ll}
 c_1 & =3\\
 c_2 & =5\\
 c_3 & =-3\\
 c_4 & =0.5\\
 c_5 & =-1\\
 c_6 & =2\\
 c_7 & =-0.7\\
 c_8 & =5.5\\
 c_9 & =-3\\
 c_{10} & =2
 \end{array}
\]
One could choose the functions $k_i(t)$ to be equal to some
constants over $[0,\infty)$. This however would require knowledge
of the exact value for a width of the dead-zone (parameter
$\delta$) in the adjustment algorithm for this particular set of
$k_i(t)$.

We simulated tracking system (\ref{track_ex}) with algorithm
(\ref{alg_ex}) for 400 trials, choosing the initial conditions for
the estimates $\hat{\alpha}(0)$ randomly in the hypercube
$[0,12]^{10}$ for every trial, initial conditions for $k_i(t)$
were set to zero. Each trial consisted of $10000$ periods (epoch)
and each epoch lasted for $T+\Delta T_2=3$ seconds. In order to
check the sensitivity of the approach to the numerical integration
we used a simple Euler's method of the first order with
integration step $\delta t=0.0001$ seconds to approximate the
solutions of $\hat{x}_i(t)$, $\hat{\alpha}_i(t)$ and $k_i(t)$. In
order to judge effectiveness of our algorithm we introduced the
following criteria:
\[
d(t)=\sqrt{\sum_{i=1}^{10}(\hat{\alpha}_i(t)-{\alpha}_i)^{2}}
\]
\[
R(t)=\sum_{i=0}^{(T+\Delta T_2)/\Delta t}\frac{e(t-T-\Delta
T_2+i\Delta t)^{2}\Delta t}{T+\Delta T_2}
\]
The histograms of distributions of distances $d(t)$ and
performance indices $R(t)$ computed in the end of each trial are
shown in Fig. 6 and 7, respectively (we made sure that
$d(0)-d((T+\Delta T_2)10000)>0$ for every trial). It can be clearly
seen from the figures  that after application of the algorithm
(\ref{alg_ex}) the distributions of the distances $d((T+\Delta
T_2)10000)$ and $R((T+\Delta T_2)10000)$ are significantly
shifted to the left towards zero.

\section{Conclusion}

In this work the problem of estimating the parameters for a
function represented by sigmoid superposition has been analyzed.
The key to our proposal is the transformation of this static
nonlinearity into a linear combination of  solutions of a system
of differential equations. These equations are linear in
parameters but nonlinear with respect to the state variables. We
considered the dynamics of an unperturbed system of differential
logistic equations. It was found that a linear combination of the
system solutions may realize any continuous function over interval
$[0,T]$ with given tolerance $\varepsilon>0$. This tolerance can
be made arbitrary small as a function of the number of equations,
with corresponding parameters and initial conditions. In addition,
we showed that a system of logistic equations with time-varying
parameters can realize a function with multiple inputs. The
results enabled us to consider a system with coupled equations via
output function $y(\hat\bfx)$ as a generator of almost any
dynamical system as long as it is smooth in its state and output
variables.

The linearity of the resulting  system with respect to its unknown
parameters allowed us to apply conventional methods and ideas of
adaptive control in order to estimate their values for a given
reference function. Extension of both the reference and tracking
signals to be repeatable (periodic) over $[0,\infty)$ interval
played a crucial role in our analysis. This feature makes it
possible to use known matching conditions (or certainty
equivalence) to design the adaptation algorithms. Stability
analysis has been performed for the learning schemes introduced.

The current algorithm is able to produce the estimates that
approach the true values of unknown system parameters within a
bounded domain. However, convergence to these true values cannot
be guaranteed. It should be mentioned, however, that the problem
of finding a flawless algorithm is all but solved by our proposal.
The most difficult hurdles to knock down were shown to be the
boundedness of solutions and the problem of determining the
maximum amplitude of unmodeled dynamics (when the reference signal
is not exactly a superposition of sigmoid function). Though we
offered possible solution to these issues in the present paper,
more effective ones may still exist. Finding these may be a topic
for future research.




\bibliographystyle{plain}
\bibliography{LogisticModel_ArXiv}

\section{Appendix}

{\it Theorem \ref{logistic_approximation} proof}. We prove the
theorem in 3 steps. First, we transform the original system
(\ref{logistic_equations}) into a system with its right-hand side
depending on one set of parameters
($\alpha=(\alpha_1,\dots,\alpha_n)^{T}$ only instead of the two
sets $\alpha$ and $\beta$). Second, for each $x_i$, $i\in
\{1,\dots,n\}$ we show that the solution $x_i(t)$ belongs to the
interval $[0,1]$ for any $x_i(0)\in(0,1)$; $x(t)$ is a  monotonic
and sigmoidal function with parameters depending on $\alpha$ and
initial conditions. Therefore, to conclude the proof it is
sufficient to apply a widely-known result\footnote{\it Let $f$ be
any continuous sigmoidal function. Then finite sums of the form: $
\sum_{i=1}^n c_i f(a_i \bfx + b_i)$, $ a_i\in R^n$, $\bfx\in R^n$,
$b_i\in R$ are dense in $C(I_n)$.} from approximation theory
\cite{Cybenko};\cite{Funahashi}.

Let us start with
\begin{lem}\label{lemma_transform}
Let system (\ref{logistic_equations}) be given and $\beta_i\neq
0$. Then there is a linear transformation $\hat{x}_i=\beta_i x_i$
of system (\ref{logistic_equations}) coordinates that the
following holds:
\begin{eqnarray}\label{logistic_equations1_1}
\dot{\hat x}_1&=&\alpha_1 \hat x_1 (1-\hat x_1); \nonumber \\
\dot{\hat x}_2&=&\alpha_2 \hat x_2 (1- \hat x_2); \nonumber \\
\cdots&=& \cdots \nonumber \\
\dot{\hat x}_n&=&\alpha_n \hat x_n (1-\hat x_n); \nonumber \\
y({\bfx})&=&\mathbf{\hat C}^T \bfx=\sum_{i}\frac{C_i}{\beta_i}
\hat{x}_i, \ \ \hat{x}_i(0)=\beta_i\Delta_i,
\end{eqnarray}
\end{lem}
{\it Lemma \ref{lemma_transform} proof}. The  proof is a routine
procedure. Let us  calculate $\dot{\hat x}_i=\beta_i \dx_i$:
\[ \dot{\hat
x}_i=\beta_i \dx_i=\alpha_i \beta_i x_i (1 - \beta_i x_i)=\alpha_i
{\hat x}_i (1-{\hat x}_i).
\]
The rest of the lemma proof is quite obvious and we skipped it.
{\it The lemma is proven.}\\
\begin{rmk}\label{trans_rem_C}\normalfont
It is desirable to note that the linear transformation
$\hat{x}_i=\beta_i x_i$ is one-to-one, and for any system
(\ref{logistic_equations1_1}) we can derive its transformed
version in the form of system (\ref{logistic_equations1}) by the
inverse transformation $x_i=1/\beta_i \hat{x}_i$. Therefore in the
rest of the proof we will deal with system
(\ref{logistic_equations1_1}). In addition, it is always possible to make a transformation such that the
resulting $\alpha_i$ will be positive. Furthermore, given system
(\ref{logistic_equations}), one can choose such linear
transformation $\hat{x}_i=1/C_i x_i$ that the transformed system
obeys
\begin{eqnarray}\label{logistic_equations1_2}
\dot{\hat x}_1&=&\alpha_1 \hat x_1 (1-\beta_1 C_1\hat x_1); \nonumber \\
\dot{\hat x}_2&=&\alpha_2 \hat x_2 (1- \beta_2 C_2 \hat x_2); \nonumber \\
\cdots&=& \cdots \nonumber \\
\dot{\hat x}_n&=&\alpha_n \hat x_n (1-\beta_n C_n \hat x_n); \nonumber \\
y({\bfx})&=&\mathbf{C}^T\bfx = \sum_{i}\frac{C_i}{C_i}
\hat{x}_i=\sum_{i}\hat{x}_i, \ \ \hat{x}_i(0)=\Delta_i/C_i,
\end{eqnarray}
thus eliminating the parametric uncertainties in output function
$y(\hat{\bfx})$ and replacing them by the parametric uncertainties
of linearly parameterized system (\ref{logistic_equations1_2})
with known output function $y(\hat\bfx)$.
\end{rmk}

 Let us consider the properties of
each $i$-th equation of system (\ref{logistic_equations1_1}). We
formulate the next lemma:

\begin{lem}\label{lemma1} Let the following differential equation be given:
\begin{equation}\label{log_1}
\dx=k x (1-x), \ \ k\neq 0,
\end{equation}
and $x(t)$ is a solution of system (\ref{log_1}) for initial
condition $x(0)=x_0$, $x_0\in (0,1)$. Then the next statements hold for equation (\ref{log_1}): \\
1) $x(t)$ is a monotonic function with respect to $t>0$; \\
2) $x(t)\rightarrow 1$ at $t\rightarrow\infty$ for $k>0$ and
$x_0\in (0,1)$; $x(t)\rightarrow 0$ at $t\rightarrow\infty$
for $k<0$ and $x_0\in (0,1)$ \\
3) $x(t)$ is unique for any $t>0$ and initial condition $x_0\in
(0,1)$.
\end{lem}
{\it Lemma \ref{lemma1} proof}. Statement 1) of the lemma proof is
obvious and therefore has been skipped here (see, for example
\cite{Strogatz}). Let us prove statement 2) of the lemma.  We consider the following function:
\begin{equation}\label{l1_proof}
V(x)=0.5(x-1)^2.
\end{equation}
It is clear that function $V(x)$ is well-defined and positive
definite for any $x>0$. Moreover, $V(x)\rightarrow\infty$ at
$x\rightarrow\infty$ and $V(x)=0$ at $x=1$. These facts allow us
to consider function $V$ as Lyapunov's candidate for system
(\ref{log_1}). Let us calculate $\dot{V}$:
\[
\dot{V}=(x-1)\dx = -kx(1-x)^2\leq 0.
\]
We observe that $V>0$ and $\dot{V}=-kx(1-x)^2<0$ for $x>0$, $x\neq
1$. For any $x\in(0,1)$, $V(x(0))-V(x(t))>0$ and therefore
$x(t)>x(0)$. Hence the next inequality holds:
\[
\dot{V}=(x-1)\dx \leq -kx(0)(1-x)^2.
\]
This can be written as follows:
\[
\dot{V} \leq-kx(0)2V(x).
\]
Hence $V\rightarrow 0$ asymptotically, and $x(t)\rightarrow 1$ at
$t\rightarrow\infty$ for any $x\in (0,1)$. To prove the second
part of statement 2, where $k<0$, it is sufficient to consider the
following Lyapunov's candidate $V(x)=0.5x^2$. Its derivative
satisfies the following equation: $\dot{V}(x)=kx^2(1-x)$ and is
obviously negative definite over $x\in [0,1)$.

 Uniqueness of $x(t)$
follows directly from the continuity of equation (\ref{log_1})
right part \cite{Pontryagin}. {\it Lemma \ref{lemma1} is
proven}.\\

Regarding lemma \ref{lemma1}, we observe that
system (\ref{logistic_equations1_1}) solutions for $\alpha_i>0$
are completely defined by the choice of initial conditions
$\hat{x}_i(0)$. This means that if $\hat{x}_i(t+\tau)$ and
$\breve{x}_i(t)$ are solutions of system
(\ref{logistic_equations1_1}) and
$\hat{x}_i(t+\tau)=\breve{x}_i(t)$ for any $t\geq 0$, then
\[
\hat{x}_i(t+\tau)=\breve{x}_i(t) \Leftrightarrow
\hat{x}_i(\tau)=\breve{x}_i(0).
\]
In other words, for each solution $\hat{x}_i(t)$ time-shift is
equivalent to choice of initial conditions. Moreover, it is easy
to see that for any $\tau \in (-\infty,\infty)$ and $\hat{x}_i(0)
\in (0,1)$ there is an initial condition $\breve{x}_i(0)$ such
that $\hat{x}_i(t+\tau)=\breve{x}_i(t)$.

All we have to prove now is that $\hat{x}_i(t)$ is a sigmoidal
function. Let us consider $\dot{\hat x}_i$. As it follows from
system (\ref{logistic_equations1_1}) equations, $\hat{x}_i(t)$
time-derivative is:
\[
\frac{\pd \hat{x}_i(t)}{\pd t}=\alpha_i {\hat x}_i(t)(1-{\hat
x}_i(t)),
\]
then
\begin{equation}\label{sigm}
\hat{x}_i(t)=\int \alpha_i {\hat x}_i(t)(1-{\hat x}_i(t)) dt
=f(\alpha_i t + b_i) + D,
\end{equation}
where
\[
f(\alpha_i t + b_i)=\frac{1}{1+e^{-(\alpha_i t + b_i)}}, \  D=0.
\]
As initial conditions of system (\ref{logistic_equations1_1})
completely define time-shifts of the solutions  $\hat{x}_i(t)$,
 coefficients $b_i$ in (\ref{sigm}) depend on initial
conditions $\hat{x}_i(0)$ only.

We just proved that $i$-th solution of system
(\ref{logistic_equations1_1}) can be written in the following
manner:
\[
\hat{x}_i(t)= f( \alpha_i t +b_i),
\]
where $b_i\in(-\infty,\infty)$, $b_i=f^{-1}(\hat{x}_i(0))$ depends
on $\hat{x}_i(0)\in (0,1)$ explicitly and $f(\cdot)$ is the
sigmoid function. Let us consider output $y(\bfx)$ of system
(\ref{logistic_equations1_1}):
\[
y(\bfx)=\sum_{i=1}^n \left(\frac{C_i}{\beta_i} f( \alpha_i t
+b_i)\right).
\]
We denote $\hat{c}_i={C_i}/\beta_i$, so $y(\bfx)$ can be written
in the form:
\[
\sum_{i=1}^n \left(\hat{c}_i f(\alpha_i t +b_i)\right).
\]
Therefore, due to \cite{Cybenko}, for any $\varepsilon>0$ and
$g(t)\in C^1_{[0,T]}$ there are such $n$, $\hat{c}_i$ and $b_i$
that the following inequality holds:
\[
|\sum_{i=1}^n \left(\hat{c}_i f(\alpha_i t
+b_i)\right)-g(t)|\leq\varepsilon
\]
for $t\in [0,T]$. To conclude the proof, it is sufficient to
notice that parameters $\alpha_i$, $\beta_i$ and initial
conditions $\Delta_i$ can be restored from $b_i$ and $\hat{c}_i$.
{\it The theorem is proven.}

{\it Lemma \ref{CEquiv} proof.} The lemma proof is trivial.
Trajectories $\bfx(t)$ and $\hat\bfx(t)$ of (\ref{ONO2}) are
bounded, then sum
\[
|\hat\mathbf{C}^T\sum_{i=1}^{n}\left(\alpha_i(\xi_{1,i}(\hat{\bfx})-\xi_{1,i}(\bfx))+\beta_i(\xi_{2,i}(\hat{\bfx})-\xi_{2,i}(\bfx))\right)(1-\lambda(t,D))|<D_2,
\]
where $D_2>0$. Therefore the coefficients
$k_i^{\ast}$ (if exist) should satisfy the following inequality
\[
\frac{D_2}{\varepsilon}<\frac{D_2}{\delta} <-\sum_{i=1}^n k_i^\ast
\hat{c}_i
\]
for $\varepsilon>\delta>0$.  Vector $\hat{\mathbf{C}}\neq 0$, hence
there exists at least one $\hat{c}_i\neq 0$. Therefore
there exists at least one  vector
$k^{\ast}=(k^{\ast}_1,\dots,k^{\ast}_n)^{T}$ such that
\[
\hat{\mathbf{C}}^T k^{\ast}<-\frac{D_2}{\delta}
\]
Therefore inequality (\ref{CEquiv_inequality}) is satisfied for
every $\epsilon>\delta>0$. {\it The lemma is proven.}

{\it Theorem \ref{learning_theorem} proof.} According to the
theorem assumptions vector $\hat{\mathbf{C}}\neq 0$. Therefore, from
Lemma \ref{CEquiv} it follows that there exist coefficients
$k^{\ast}_i$ such that
\begin{eqnarray}\label{lt_1}
|\hat\mathbf{C}^T\sum_{i=1}^{n}\left(\alpha_i(\xi_{1,i}(\hat{\bfx})-\xi_{1,i}(\bfx))+\beta_i(\xi_{2,i}(\hat{\bfx})-\xi_{2,i}(\bfx))\right)(1-\lambda(t,D))|+\epsilon\sum_{i=1}^n
k_i^{\ast}\hat{c}_i<0
\end{eqnarray}
for any $\epsilon>\delta-\delta_1$, where $\delta>\delta_1>0$.
Define the following set of time intervals:
\begin{eqnarray}
\Delta_{t,0}&=&\{\Delta
({2i,2i+1})=[t_{2i},t_{2i+1}]|\lambda(t,D)=0 \ \forall
t\in[t_{2i},t_{2i+1}],  \nonumber \\ & & i\in\mathcal{N}, \
t_0<t_1\dots<t_{j}<t_{j+1}<t_{j+2}<\dots \}.\nonumber
\end{eqnarray}
\begin{eqnarray}
\Delta_{t,1}&=&\{\Omega
({2i+1,2i+2})=(t_{2i+1},t_{2i+2})|\lambda(t,D)=1 \ \forall
t\in(t_{2i+1},t_{2i+2}), \nonumber \\
& & t_1<t_2\dots<t_{j}<t_{j+1}<t_{j+2}<\dots \}.\nonumber
\end{eqnarray}

Consider the following positive-definite function
\[
V(e,\hat\alpha,\hat{\beta},K)=\int_{0}^{e}S_{\delta}(\nu) \nu d\nu
+ 0.5\|\hat{\alpha}-\alpha\|^{2}_{\gamma^{-1}} +
0.5\|\hat{\beta}-\beta\|^{2}_{\gamma^{-1}} +
0.5\|K(t)-k^{\ast}\|^{2}_{\gamma^{-1}},
\]
where $k^{\ast}_i$ satisfy inequality (\ref{lt_1}) for every
$\epsilon>\delta-\delta_1$. Its time-derivative over the set
$\Delta_{t,0}$ can be expressed as follows
\[
\frac{d}{dt}V(e,\hat\alpha,\hat{\beta},K)=S_{\delta}(e)\left(e\dot{e}-\sum_{i=1}^{n}\left((\hat{\alpha}_i-{\alpha}_i)e\hat\mathbf{C}^T\xi_{1,i}(\hat\bfx)-
(\hat{\beta}_i-{\beta}_i)e\hat\mathbf{C}^T\xi_{2,i}(\hat\bfx)-(k_i(t)-k_i^{\ast})e^2\hat{c}_i\right)\right)
\]
It is clear that $\dot{V}=0$ for any $|e|<\delta$ as
$S_{\delta}(e)\equiv 0$ for all $|e|<\delta$. Let $|e|\geq
\delta$, then
\begin{eqnarray}
\dot{V}&=&S_{\delta}(e)e\left(\hat\mathbf{C}^{T}
\sum_{i=1}^{n}\left(\hat{\alpha}_i\xi_{1,i}(\hat\bfx)-{\alpha}_i\xi_{1,i}(\bfx)+\hat{\beta}_i\xi_{2,i}(\hat\bfx)-{\beta}_i\xi_{2,i}(\bfx)\right)
+ \hat\mathbf{C}^{T}K(t)e
\right)-\nonumber\\
 & & S_{\delta}(e) e \left(\sum_{i=1}^{n}(\hat{\alpha}_i-{\alpha}_i)\hat\mathbf{C}^T\xi_{1,i}(\hat\bfx)-
(\hat{\beta}_i-{\beta}_i)\hat\mathbf{C}^T\xi_{2,i}(\hat\bfx)-(k_i(t)-k_i^{\ast})e\hat{c}_i\right)\nonumber
\\
& = &
S_{\delta}(e)e\left(\sum_{i=1}^n\hat{\mathbf{C}}^{T}\alpha_i(\xi_{1,i}(\hat\bfx)-\xi_{1,i}(\bfx))+\hat{\mathbf{C}}^{T}\beta_i(\xi_{2,i}(\hat\bfx)-\xi_{2,i}(\bfx))+k^{\ast}_i\hat{c}_i
e \right)\nonumber\\
&\leq  &
S_{\delta}(e)|e|\left(|\left(\sum_{i=1}^n\hat{\mathbf{C}}^{T}\alpha_i(\xi_{1,i}(\hat\bfx)-\xi_{1,i}(\bfx))+\hat{\mathbf{C}}^{T}\beta_i(\xi_{2,i}(\hat\bfx)-\xi_{2,i}(\bfx))\right)|+\sum_{i=1}^nk_i^{\ast}\hat{c}_i|e|\right)\nonumber\\
& \leq
&S_{\delta}(e)|e|\sum_{i=1}^{n}{k^{\ast}_i}\hat{c}_i(|e|-\delta+\delta_1)\leq
S_{\delta}(e)|e|\sum_{i=1}^{n}{k^{\ast}_i}\hat{c}_i\delta_1\leq 0.
\end{eqnarray}
(In order to get the last inequality note that sum
$\sum_{i=1}^nk^{\ast}_i\hat{c}_i$ must be negative.) Taking into
account that $\dot{V}$ is not positive over $[t_{2i},t_{2i+1}]$
and that
$e(t_i)=0$ (because the states of both reference and tracking
systems are forced to move to $\bfx(0)$ over $\Delta_{t,1}$), one
can write
\begin{eqnarray}
&  &
V(e(t_{2i}),\hat{\alpha}(t_{2i}),\hat{\beta}(t_{2i}),\hat{K}(t_{2i}))-V(e(t_{2i+1}),\hat{\alpha}(t_{2i+1}),\hat{\beta}(t_{2i+1}),\hat{K}(t_{2i+1}))\nonumber
\\
& = &  0.5\|\hat{\alpha}(t_{2i})-\alpha\|^{2}_{\gamma^{-1}} +
0.5\|\hat{\beta}(t_{2i})-\beta\|^{2}_{\gamma^{-1}} +
0.5\|K(t_{2i})-k^{\ast}\|^{2}_{\gamma^{-1}} - \nonumber \\
& & 0.5\|\hat{\alpha}(t_{2i+1})-\alpha\|^{2}_{\gamma^{-1}} -
0.5\|\hat{\beta}(t_{2i+1})-\beta\|^{2}_{\gamma^{-1}} -
0.5\|K(t_{2i+1})-k^{\ast}\|^{2}_{\gamma^{-1}}\nonumber\\
& > &
\int_{t_{2i}}^{t_{2i+1}}S_{\delta}(e)|e(\tau)\sum_{j=1}^{n}{k^{\ast}_j}\hat{c}_j|\delta_1
d\tau + \int_0^{e(t_{2i+1})}S_{\delta}(\nu)\nu d\nu. \nonumber
\end{eqnarray}

Consider the following series:
\begin{eqnarray}
W(n) &=& 0. 5 \sum_{i=0}^{n}\left(
\|\hat{\alpha}(t_{2i})-\alpha\|^{2}_{\gamma^{-1}} +
\|\hat{\beta}(t_{2i})-\beta\|^{2}_{\gamma^{-1}} +
\|K(t_{2i})-k^{\ast}\|^{2}_{\gamma^{-1}} -\right. \nonumber \\
& & \left. \|\hat{\alpha}(t_{2i+1})-\alpha\|^{2}_{\gamma^{-1}} -
\|\hat{\beta}(t_{2i+1})-\beta\|^{2}_{\gamma^{-1}} -
\|K(t_{2i+1})-k^{\ast}\|^{2}_{\gamma^{-1}}\right).\nonumber
\end{eqnarray}
One can notice that
\begin{eqnarray}
& & \|\hat{\alpha}(t_{2i+1})-\alpha\|^{2}_{\gamma^{-1}} +
\|\hat{\beta}(t_{2i+1})-\beta\|^{2}_{\gamma^{-1}} +
\|K(t_{2i+1})-k^{\ast}\|^{2}_{\gamma^{-1}}\nonumber \\
& = & \|\hat{\alpha}(t_{2i+2})-\alpha\|^{2}_{\gamma^{-1}} +
\|\hat{\beta}(t_{2i+2})-\beta\|^{2}_{\gamma^{-1}} +
\|K(t_{2i+2})-k^{\ast}\|^{2}_{\gamma^{-1}}\nonumber
\end{eqnarray}
as vectors $\hat{\alpha}$, $\hat{\beta}$ and $K$ remain constant
over intervals $\Delta_{t,1}$.  Therefore
\begin{eqnarray}\label{lt_2}
W(n)&=& 0.5\left(\|\hat{\alpha}(t_0)-\alpha\|^{2}_{\gamma^{-1}} +
\|\hat{\beta}(t_0)-\beta\|^{2}_{\gamma^{-1}} +
\|K(t_0)-k^{\ast}\|^{2}_{\gamma^{-1}} -\right. \nonumber \\
& & \left.\|\hat{\alpha}(t_{2n+1})-\alpha\|^{2}_{\gamma^{-1}} -
\|\hat{\beta}(t_{2n+1})-\beta\|^{2}_{\gamma^{-1}} -
\|K(t_{2n+1})-k^{\ast}\|^{2}_{\gamma^{-1}}\right)\nonumber\\
&>&\sum_{i=0}^n\int_{t_{2i}}^{t_{2i+1}}S_{\delta}(e)|e(\tau)\sum_{j=1}^{n}{k^{\ast}_j}\hat{c}_j|\delta_1
d\tau + \sum_{i=0}^n\int_0^{e(t_{2i+1})}S_{\delta}(\nu)\nu d\nu >
0.
\end{eqnarray}
Given that  $\bfx(t)$, $\hat{\bfx}(t)$ are bounded we can conclude
that $\dot{\hat\alpha}$, $\dot{\hat\beta}$ and $\dot{K}(t)$ are
bounded and hence $\hat\alpha$, $\hat{\beta}$, $K(t)$ are bounded.
Furthermore, the following inequality holds
\begin{eqnarray}
& & 0.5\left(\|\hat{\alpha}(t_0)-\alpha\|^{2}_{\gamma^{-1}} +
\|\hat{\beta}(t_0)-\beta\|^{2}_{\gamma^{-1}} +
\|K(t_0)-k^{\ast}\|^{2}_{\gamma^{-1}}\right)
>\sum_{i=0}^n\int_{t_{2i}}^{t_{2i+1}}S_{\delta}(e)|e(\tau)\sum_{j=1}^{n}{k^{\ast}_j}\hat{c}_j|\delta_1
d\tau \nonumber \\
& = &
\int_0^{t_{2n+1}}S_{\delta}(e)\lambda(\tau,D)|e(\tau)\sum_{j=1}^{n}{k^{\ast}_j}\hat{c}_j|\delta_1d\tau>0.\nonumber
\end{eqnarray}
Hence
\[
0<\int_0^{\infty}S_{\delta}(e)\lambda(\tau,D)|e(\tau)\sum_{j=1}^{n}{k^{\ast}_j}\hat{c}_j|\delta_1d\tau<\infty.
\]
Let us consider the following time-intervals
$\Delta_i=[\tau_{2i},\tau_{2i+1}]: \ |e|\lambda(t,D)\geq\delta  \
\forall t\in \Delta_i$, $i\in\{0,1,\dots,\infty\}$. As
$|e(t)|>\delta$ it is clear that
\begin{eqnarray}
\infty&>&\int_0^{\infty}S_{\delta}(e)\lambda(\tau,D)|e(\tau)\sum_{j=1}^{n}{k^{\ast}_j}\hat{c}_j|\delta_1d\tau>
\int_0^{\infty}S_{\delta}(e)\lambda(\tau,D)|\delta
\sum_{j=1}^{n}{k^{\ast}_j}\hat{c}_j|\delta_1d\tau\nonumber \\
& =& \sum_{i=0}^{\infty}\Delta_i |\delta
\sum_{j=1}^{n}{k^{\ast}_j}\hat{c}_j|\delta_1 >0. \nonumber
\end{eqnarray}
Then series
\[
\sum_{i=0}^{\infty}\Delta_i
|\delta\sum_{j=1}^{n}{k^{\ast}_j}\hat{c}_j|\delta_1
\]
converges and, therefore, $\Delta_i\rightarrow 0$ as
$i\rightarrow\infty$. In order to finish the proof of the theorem,
it is sufficient to consider the error function $e(t)$ over
intervals $\Delta_i$. Derivative $\dot{e}$ is bounded (say
$|\dot{e}|<D_3$) as vectors $\bfx$, $\hat{\bfx}$, $\hat{\alpha}$,
$\hat{\beta}$, $K(t)$ are bounded. Therefore for any
$t\in\Delta_i$:
\begin{eqnarray}
|e(t)\lambda(t,D)|&=&|e(\tau_{2i})+\int_{\tau_{2i}}^{\tau_{2i+1}}\dot{e}(\tau)d\tau|\leq|e(\tau_{2i})|+|\int_{\tau_2i}^{\tau_{2i+1}}\dot{e}(\tau)d\tau|\nonumber\\
& \leq & |e(\tau_{2i})| + |\Delta_i|D_3=\delta+\Delta_i D_3.
\nonumber
\end{eqnarray}
Then
\[
\lim_{t\rightarrow \infty} \sup |e(t)\lambda(t,D)|=\delta.
\]
Hence for any arbitrary small $\delta_1>0$ there exists such $t_1$
that
\[
|e(t)\lambda(t,D)|<\delta+\delta_1
\]
for any $t>t_1$. {\it The theorem is proven.}

\begin{figure}\label{fig1}
\begin{center}
\includegraphics
[width=\columnwidth]{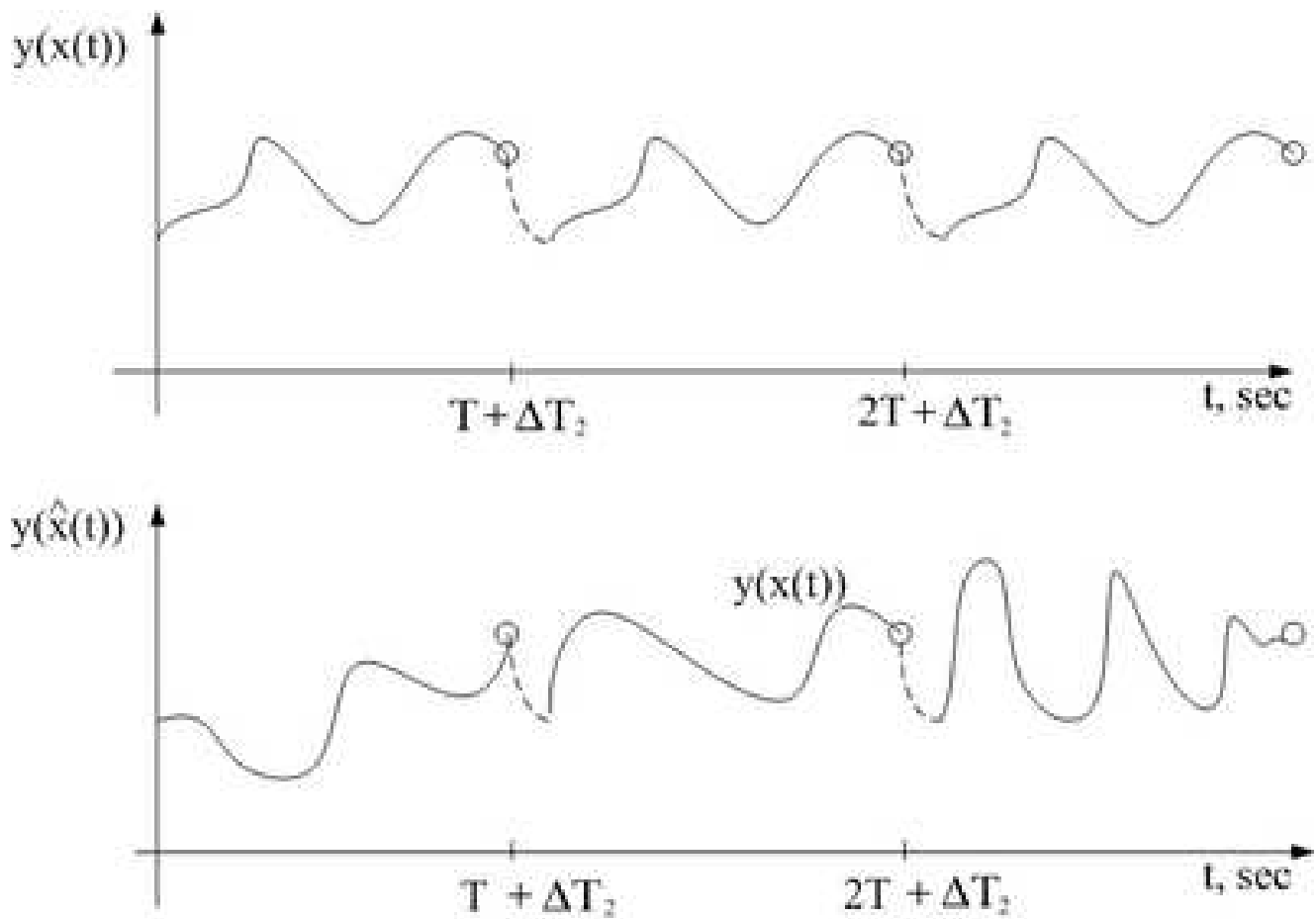} \caption{Periodical
extension of the reference signal $g(t)$ defined over $[0,T]$}
\end{center}
\begin{center}
\end{center}
\end{figure}

\begin{figure}\label{fig2}
\begin{center}
\includegraphics
[width=\columnwidth]{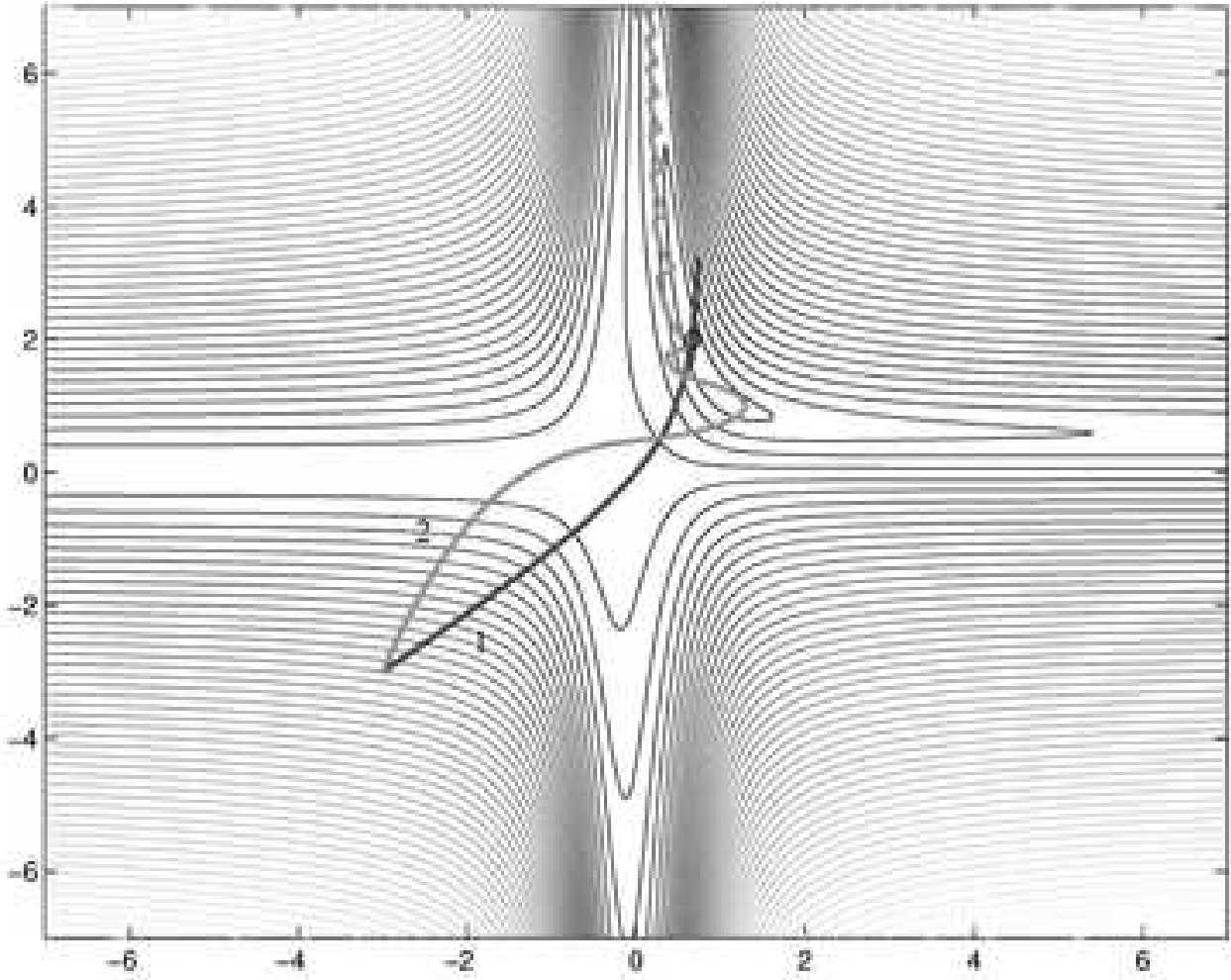} \caption{
Trajectories $\hat{\alpha}(t),\hat{c}(t)$ in system
(\ref{example3}) with algorithm (\ref{example3_learning}) (Curve
1) and algorithm (\ref{example3_gradient}) (Curve 2) starting from
point $(-3,-3)$. Global minimum is marked by circle}
\end{center}
\begin{center}
\end{center}
\end{figure}

\begin{figure}\label{fig3}
\begin{center}
\includegraphics
[width=\columnwidth]{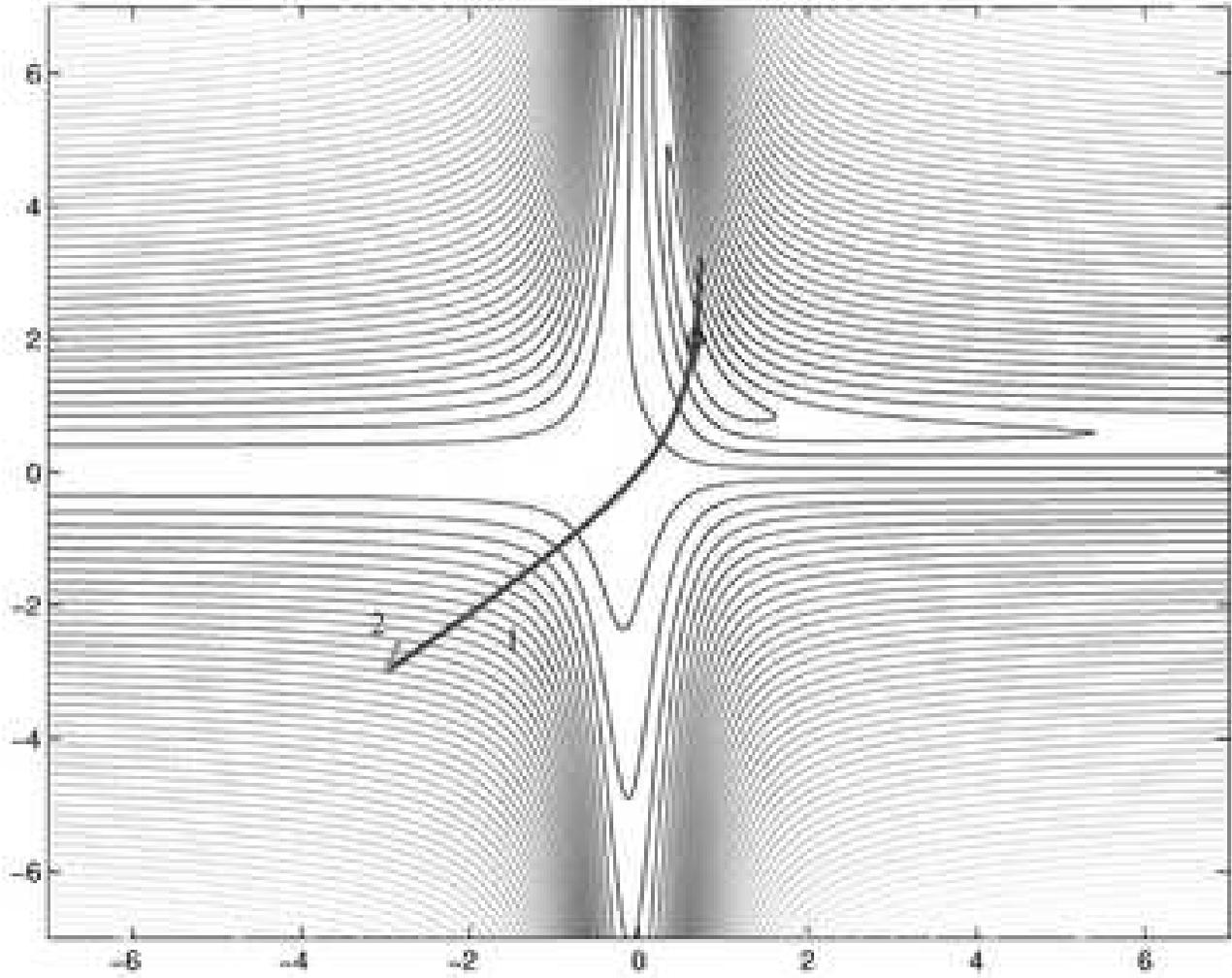} \caption{
Trajectories $\hat{\alpha}(t),\hat{c}(t)$ in system
(\ref{example3}) with algorithm (\ref{example3_learning}) (Curve
1) and algorithm (\ref{example3_gradient}) (Curve 2) starting from
point $(-3,-3)$. The trajectories have been shown for time
interval $[0,900]$ sec. }
\end{center}
\begin{center}
\end{center}
\end{figure}

\begin{figure}\label{fig4}
\begin{center}
\includegraphics
[width=\columnwidth]{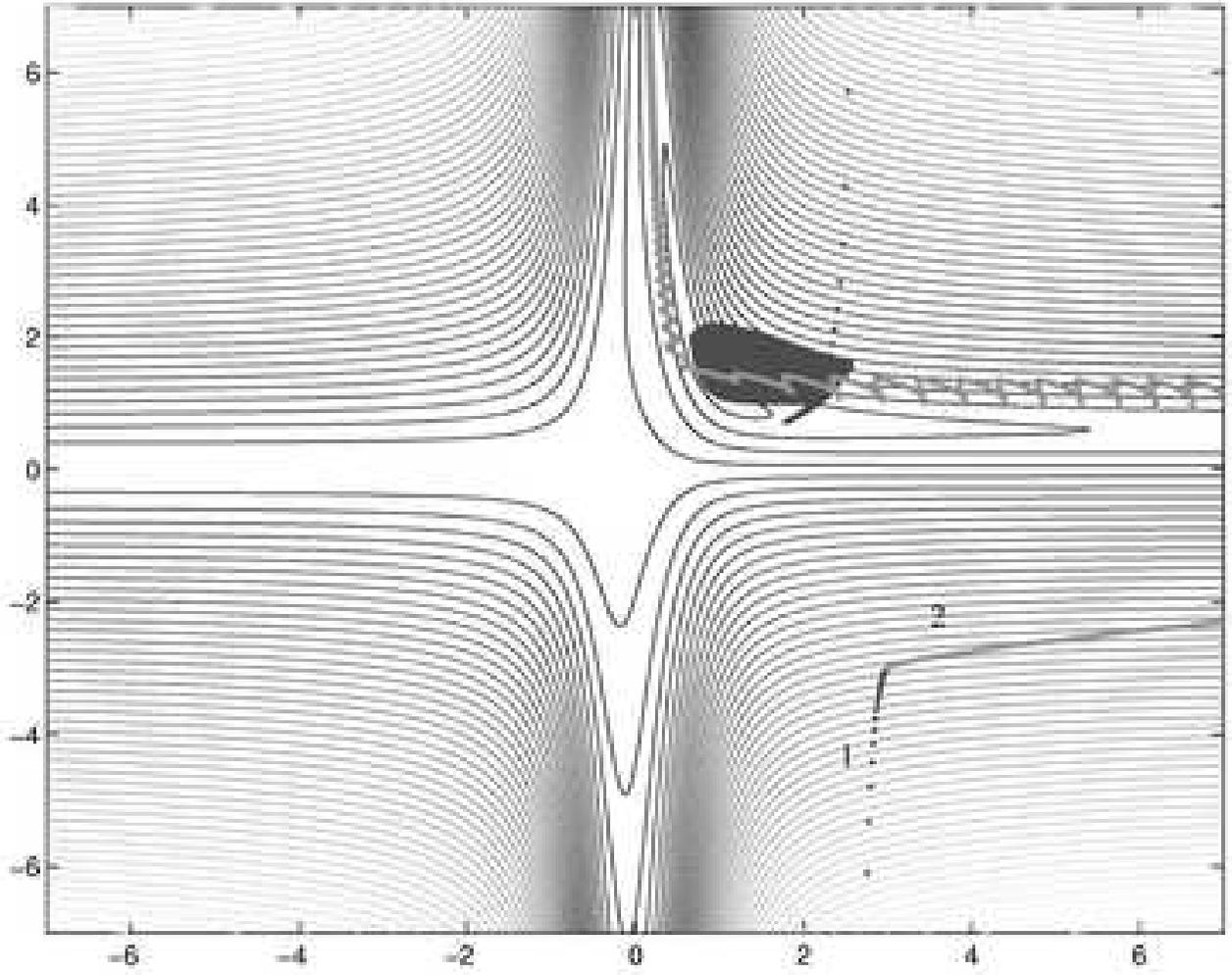} \caption{ Trajectories
$\hat{\alpha}(t),\hat{c}(t)$ in system (\ref{example3}) with
algorithm (\ref{example3_learning}) (Curve 1) and algorithm
(\ref{example3_gradient}) (Curve 2) starting from point $(3,-3)$.
Algorithm (\ref{example3_learning}) ensures that the estimates
reach a neighborhood of the global minimum in very short time and
then to approach it with oscillations in the parameter space
(blob-like part of the trajectory).}
\end{center}
\begin{center}
\end{center}
\end{figure}

\begin{figure}\label{fig5}
\begin{center}
\includegraphics
[width=\columnwidth]{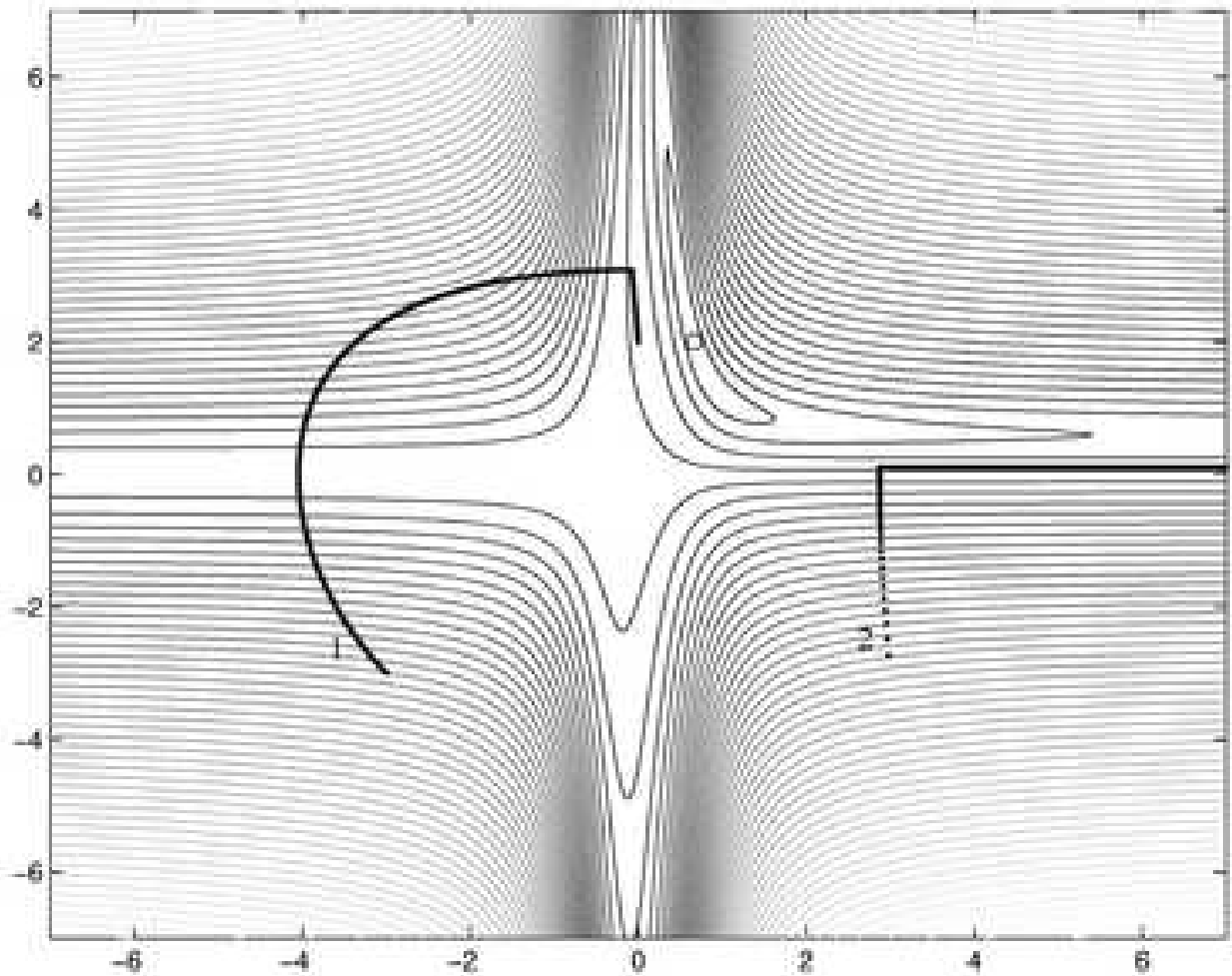} \caption{ Trajectories
$\hat{\alpha}(t),\hat{c}(t)$ in system (\ref{example3}) with batch
gradient algorithm (\ref{example3_batch}) starting from point
$(-3,-3)$ (Curve 1) and $(3,-3)$ (Curve 2). None reaches the
global minimum (marked by circle)}
\end{center}
\begin{center}
\end{center}
\end{figure}

\begin{figure}\label{fig6}
\begin{center}
\includegraphics
[width=\columnwidth]{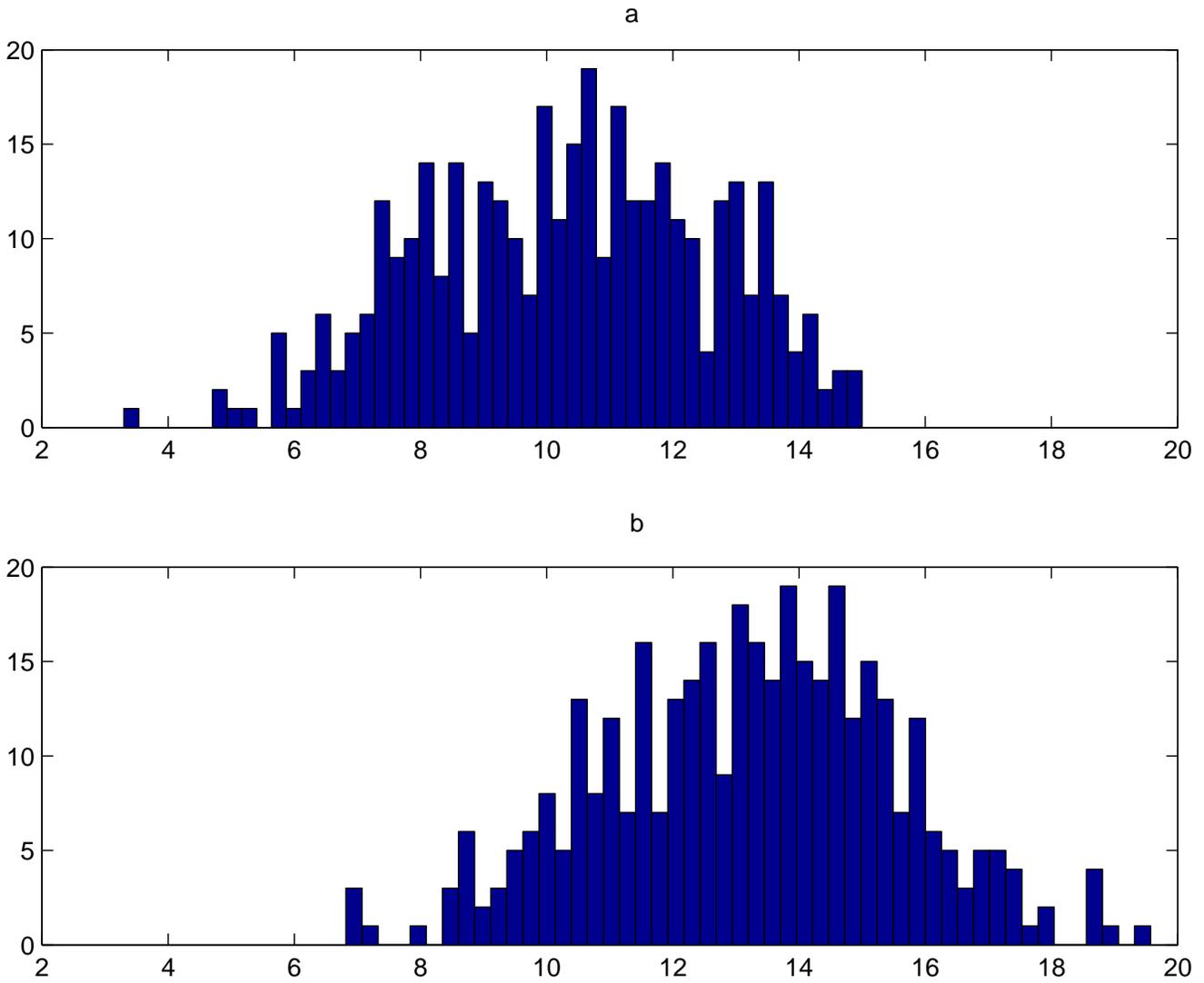} \caption{ Histograms of the
distributions of the distances $d((T+\Delta T_2)10000)$ (plot $a$)
and $d(0)$ (plot $b$) for $400$ trials with random initial
conditions for the estimates $\hat{\alpha}_i(0)$.}
\end{center}
\begin{center}
\end{center}
\end{figure}

\begin{figure}\label{fig7}
\begin{center}
\includegraphics
[width=\columnwidth]{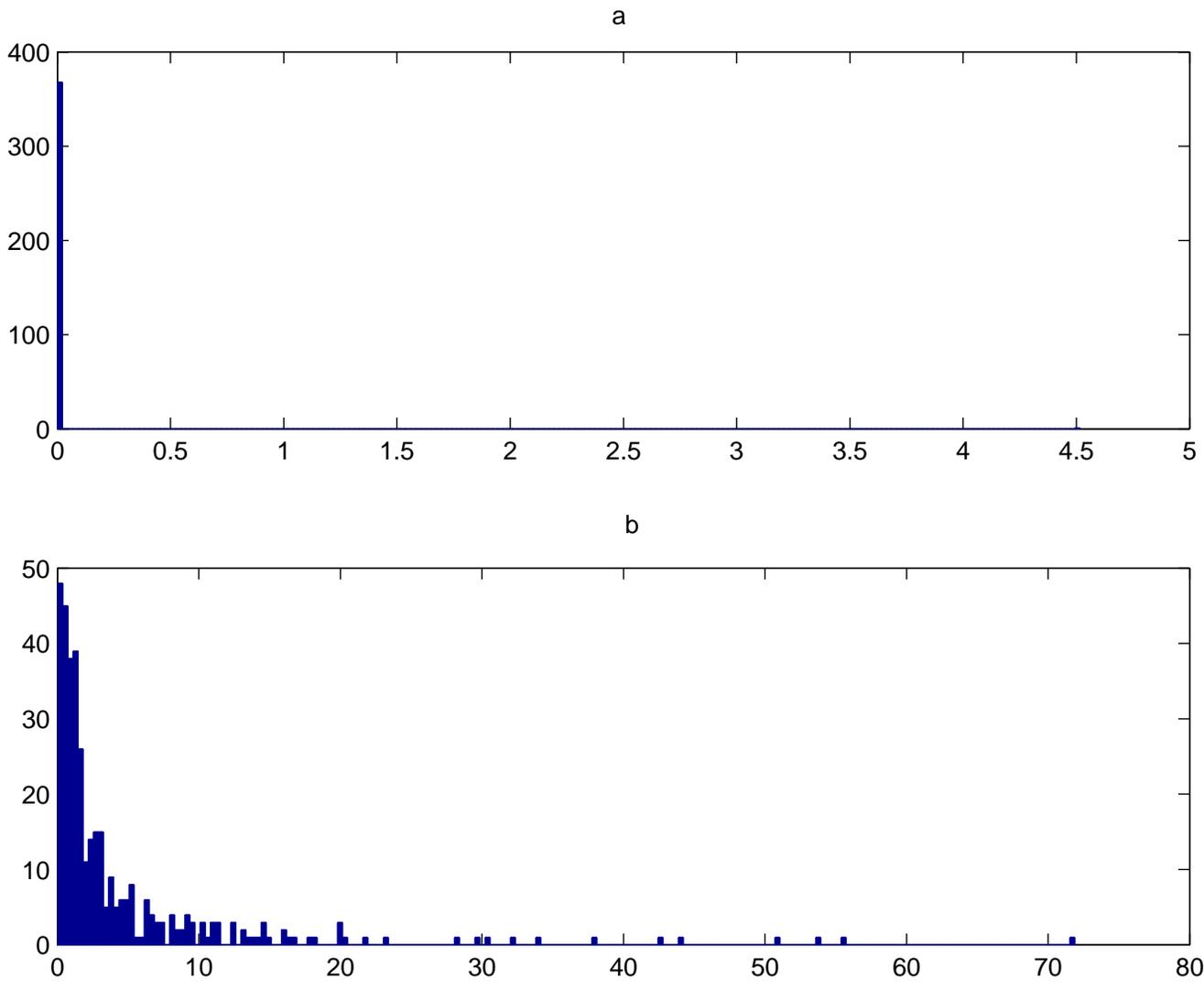} \caption{ Histograms of the
distributions of the performance indices  $R((T+\Delta T_2)10000)$
(plot $a$) and $R(0)$ (plot $b$) for $400$ trials with random
initial conditions for the estimates $\hat{\alpha}_i(0)$.}
\end{center}
\begin{center}
\end{center}
\end{figure}

\end{document}